\documentclass{gtart_a}
\pdfoutput=1
\usepackage{graphicx}
\usepackage[percent]{overpic}

\title[Ropelength bounds from quadrisecants]
  {Quadrisecants give new lower bounds\\for the ropelength of a knot}

\author{Elizabeth Denne}
\givenname{Elizabeth}
\surname{Denne}
\address{Department of Mathematics\\
Harvard University\\\newline
Cambridge\\
Massachusetts 02138\\USA}
\email{denne@math.harvard.edu}
\urladdr{}

\author{Yuanan Diao}
\givenname{Yuanan}
\surname{Diao}
\address{Department of Mathematics\\
University of North Carolina\\\newline
Charlotte\\
North Carolina 28223\\USA}
\email{ydiao@uncc.edu}
\urladdr{}

\author{John M Sullivan}
\givenname{John M}
\surname{Sullivan}
\address{Institut f\"ur Mathematik\\
MA 3--2\\Technische Universit\"at Berlin\\\newline
D--10623 Berlin\\Germany}
\email{jms@isama.org}
\urladdr{}

\doi{}

\volumenumber{10}
\issuenumber{}
\publicationyear{2006}
\papernumber{1}
\lognumber{0660}
\startpage{1}
\endpage{26}


\received{27 September 2005}
\revised{27 December 2005}
\accepted{3 January 2006}
\published{25 February 2006}
\publishedonline{25 February 2006}
\proposed{Joan Birman}
\seconded{Dave Gabai, Walter Neumann}
\corresponding{}
\editor{}
\version{}

\begin{asciiabstract}
Using the existence of a special quadrisecant line,
we show the ropelength of any nontrivial knot is at least 15.66.
This improves the previously known lower bound of 12.
Numerical experiments have found a trefoil with ropelength less
than 16.372, so our new bounds are quite sharp.
\end{asciiabstract}

\AtBeginDocument{\def\notin{\not\in}}

\newtheorem{theorem}{Theorem}[section]
 \unnumbered{maintheorem}
\newtheorem{lemma}{Lemma}[section]
\newtheorem{proposition}{Proposition}[section]
\newtheorem{corollary}{Corollary}[section]
\makeatletter
\let\c@lemma\c@theorem
\let\c@proposition\c@theorem
\let\c@corollary\c@theorem
\makeatother
\theoremstyle{definition}

\newtheorem{definition}{Definition} \unnumbered{definition}
 \unnumbered{example}
 \unnumbered{conjecture}
\newtheorem{remark}[theorem]{Remark}

\def\secn#1{\hbox{\fullref{sec:#1}}}
\def\thm#1{\hbox{\fullref{thm:#1}}}
\def\lem#1{\hbox{\fullref{lem:#1}}}
\def\cor#1{\hbox{\fullref{cor:#1}}}
\def\prop#1{\hbox{\fullref{pr:#1}}}
\def\figr#1{\hbox{\fullref{fig:#1}}}

\graphicspath{{figs/}}

\newcommand{\arc}[1]{\gamma_{#1}}
\newcommand{\len}[1]{\ell_{#1}}
\newcommand{\bdry}{\partial}

\newcommand{\setm}{\smallsetminus}
\newcommand{\eps}{\varepsilon}

\makeatletter
\newbox\overbox
\def\fakeover#1{\setbox\overbox\hbox{$#1$}\hbox
                         {$\overline{#1\hskip-\wd\overbox}$\hskip\wd\overbox}}
\def\overnoarrow#1{\mskip1.5mu\overline{\mskip-1.5mu#1}}
\def\overrightarrow#1{\mskip2mu\vbox{\m@th\ialign{##\crcr
   \rightarrowfill\crcr
   \noalign{\kern-.4pt               
        \kern-\fontdimen22\textfont2 
        \nointerlineskip}
   ${\mskip0mu\hfil\fakeover{#1}\hfil\mskip6mu}$\crcr}}\mskip-2mu}
\def\overleftrightarrow#1{\mskip-3mu\vbox{\m@th\ialign{##\crcr
   \leftarrowfill\hskip-.6em\rightarrowfill\crcr
   \noalign{\kern-.4pt               
        \kern-\fontdimen22\textfont2 
        \nointerlineskip}
   ${\mskip3mu\hfil\fakeover{#1}\hfil\mskip6mu}$\crcr}}\mskip-2mu}

\def\ray#1{{\smash{\overrightarrow{#1}}}}

\def\segment#1{{\smash{\overnoarrow{#1}}}}
\makeatother

\newcommand{\altval}{15.66}
\newcommand{\minval}{13.936} 
\newcommand{\numval}{16.372} 

\begin{document}

\begin{abstract}
Using the existence of a special quadrisecant line,
we show the ropelength of any nontrivial knot is at least~15.66.
This improves the previously known lower bound of~12.
Numerical experiments have found a trefoil with ropelength less
than~16.372, so our new bounds are quite sharp.
\end{abstract}

\maketitle

\section{Introduction}
The ropelength problem seeks to minimize the length of a knotted
curve subject to maintaining an embedded tube of fixed radius
around the curve; this is a mathematical model of tying the knot
tight in rope of fixed thickness.

More technically, the thickness~$\tau(K)$ of a space curve~$K$ is
defined by Gonzalez and Maddocks \cite{gm} to be twice the infimal radius
$r(a,b,c)$ of
circles through any three distinct points of~$K$.
It is known from the work of Cantarella, Kusner and Sullivan~\cite{cks}
that $\tau(K)=0$ unless $K$ is $C^{1,1}$, meaning that its tangent
direction is a Lipschitz function of arclength.  When $K$ is $C^1$,
we can define normal tubes around $K$, and then indeed $\tau(K)$ is the
supremal diameter of such a tube that remains embedded.  We note that in
the existing literature thickness is sometimes defined to be the radius
rather than diameter of this thick tube.

We define ropelength to be the (scale-invariant) quotient of length over
thickness.  Because this is semi-continuous even in the $C^0$ topology
on closed curves, it is not hard to show \cite{cks}
that any (tame) knot or link type has a ropelength minimizer.

Cantarella, Kusner and Sullivan gave certain lower bounds
for the ropelength of links; these are sharp in certain simple cases where
each component of the link is planar \cite{cks}.  However, these examples
are still the only known ropelength minimizers.  Recent work by
Cantarella, Fu, Kusner, Sullivan and Wrinkle \cite{cfksw}
describes a much more complicated tight (ropelength-critical)
configuration $B_0$ of the Borromean rings.  (Although the somewhat different
Gehring notion of thickness is used there, $B_0$ should still be tight, and
presumably minimizing, for the ordinary ropelength we consider here.)
Each component of $B_0$ is still planar, and it seems significantly more difficult
to describe explicitly the shape of any tight knot.

Numerical experiments by Piera\'nski \cite{pie}, Sullivan \cite{ropen} and
Rawdon \cite{raw} suggest that the minimum ropelength
for a trefoil is slightly less than~$\numval$, and that there is another tight
trefoil with different symmetry and ropelength about~$18.7$.  For comparison,
numerical simulations of the tight figure-eight knot show ropelength
just over~$21$.  The best lower bound in~\cite{cks} was $10.726$;
this was improved by Diao~\cite{diao},
who showed that any knot has ropelength more than~$12$ (meaning that
``no knot can be tied in one foot of one-inch rope'').

Here, we use the idea of quadrisecants, lines that intersect a knot
in four distinct places, to get better lower bounds for ropelength.
Almost 75 years ago, Pannwitz~\cite{pann} showed, using polygonal knots,
that a generic representative of any nontrivial knot type must have
a quadrisecant.  Kuperberg~\cite{kup}
extended this result to all knots by showing that generic knots
have essential (or topologically nontrivial) quadrisecants.
(See also the article by Morton and Mond \cite{MortonMond}.)
We will define this precisely below, as essential quadrisecants are exactly
what we need for our improved ropelength bounds.  (Note that a curve
arbitrarily close to a round circle, with ropelength thus near $\pi$, can
have a nonessential quadrisecant.)

By comparing the orderings of the four points along the knot and along
the quadrisecant, we distinguish three types of quadrisecants.
For each of these types we use geometric arguments to obtain a lower
bound for the ropelength of the knot having a quadrisecant of that type.
The worst of these three bounds is $\minval$.

In her doctoral dissertation~\cite{denne}, Denne shows that nontrivial knots
have essential quad\-ri\-se\-cants of alternating type.  This result,
combined with our \thm{alternating}, shows that any nontrivial knot has
ropelength at least $\altval$.

Nontrivial links also necessarily have quadrisecants.  We briefly consider
ropelength bounds obtained for links with different types of quadrisecants. 
These provide another interpretation of the argument showing that
the tight Hopf link has ropelength $4\pi$, as in the Gehring link problem.
But we have not found any way to improve the known ropelength estimates
for other links.

\section{Definitions and lemmas}
\begin{definition}
A \emph{knot} is an oriented simple closed curve~$K$ in~$\R^3$.
Any two points $a$ and~$b$ on a knot $K$ divide it into two complementary
subarcs $\arc{ab}$ and $\arc{ba}$.  Here~$\arc{ab}$ is the
arc from~$a$ to~$b$ following the given orientation on $K$.
If $p\in\arc{ab}$, we will sometimes write $\arc{apb}=\arc{ab}$
to emphasize the order of points along~$K$.
We will use~$\len{ab}$ to denote the length of $\arc{ab}$;
by comparison, $|a-b|$ denotes the distance from~$a$ to~$b$ in space,
the length of the segment $\segment{ab}$.
\end{definition}

\begin{definition}
An \emph{$n$--secant line} for a knot~$K$ is an oriented
line in $\R^3$ whose intersection with~$K$ has at least $n$ components.
An \emph{$n$--secant} is an ordered $n$--tuple of points in~$K$
(no two of which lie in a common straight subarc of~$K$)
which lie in order along an $n$--secant line.
We will use \emph{secant}, \emph{trisecant} and \emph{quadrisecant}
to mean $2$--secant, $3$--secant and $4$--secant, respectively.
The \emph{midsegment} of a quadrisecant $abcd$ is the segment~$\segment{bc}$.
\end{definition}

\begin{figure}[ht!]\centering
\begin{overpic}[scale=.5]{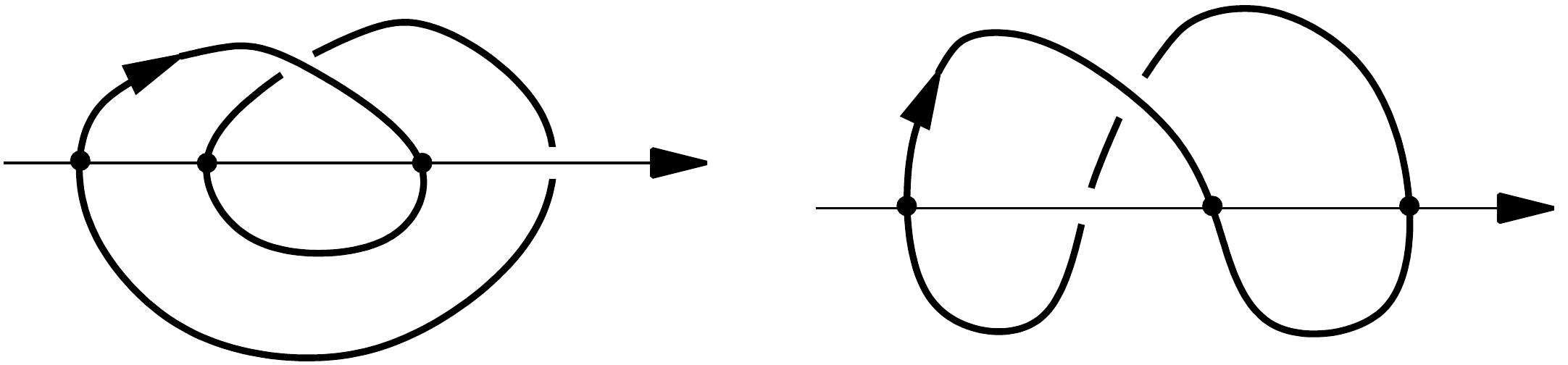}
{\small
\put(2.5,10.3){$a$}
\put(10.7,10.3){$b$}
\put(24,10.4){$c$}
\put(59.3,11){$a$}
\put(78,11){$b$}
\put(87,11){$c$}
}
\end{overpic}
\caption[Reversed and direct trisecants]{These trisecants $abc$ are reversed
(left) and direct (right) because the cyclic order of the points
along~$K$ is $acb$ and $abc$, respectively.  Flipping the orientation of
the knot or the trisecant would change its type.}
\label{fig:trisecants}
\end{figure}
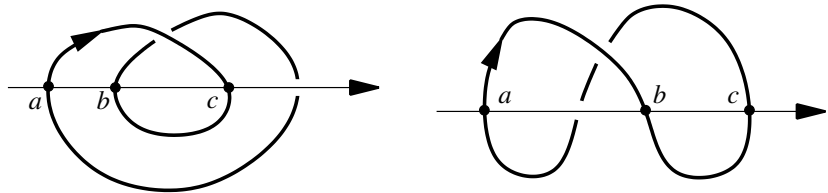

The orientation of a trisecant either agrees or disagrees
with that of the knot.  In detail, the three points of a trisecant $abc$
occur in that linear order along the trisecant line, but may occur
in either cyclic order along the oriented knot.  (Cyclic orders are cosets
of the cyclic group~$C_3$ in the symmetric group~$S_3$.)
These could be labeled by their lexicographically least elements
($abc$ and $acb$), but we choose to call them \emph{direct}
and \emph{reversed} trisecants, respectively, as in \figr{trisecants}.
Changing the orientation of either the knot or the trisecant changes
its class.  Note that $abc$ is direct if and only if $b\in\arc{ac}$.

\begin{figure}[ht!]\centering
\begin{overpic}[scale=.5]{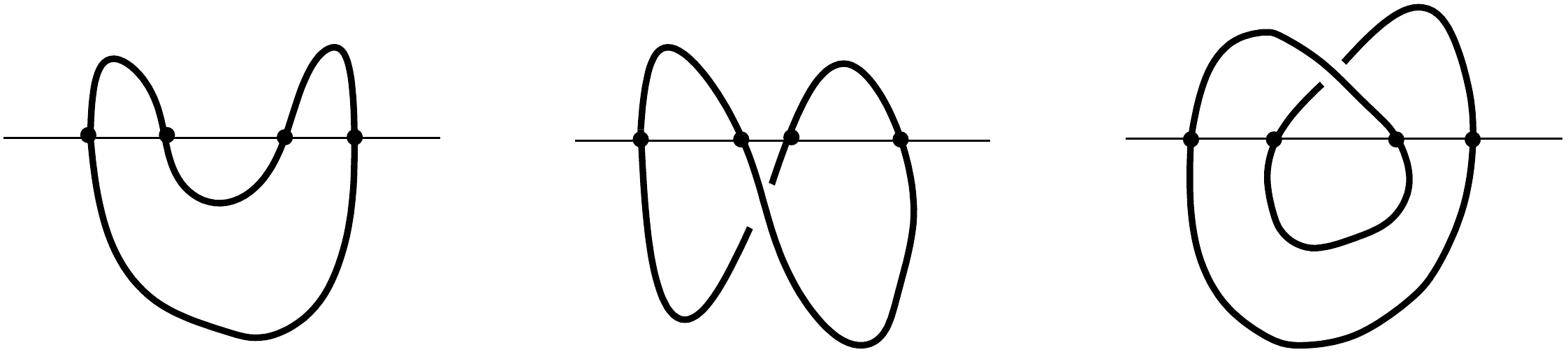}
{\small
\put(3.3,10.8){$a$}
\put(8.5,10.6){$b$}
\put(18.5,10.8){$c$}
\put(23.5,10.6){$d$}
\put(38.5,10.8){$a$}
\put(45.4,10.6){$b$}
\put(50.5,10.8){$c$}
\put(55.0,10.6){$d$}
\put(73.4,10.8){$a$}
\put(78.6,10.6){$b$}
\put(87.6,10.8){$c$}
\put(91.3,10.6){$d$}
}
\end{overpic}
\caption[Simple, flipped and alternating quadrisecants]
{Here we see quadrisecants $abcd$ on each of three knots.
From left to right, these are simple, flipped and alternating,
because the dihedral order of the points along~$K$ is $abcd$, $abdc$
and $acbd$, respectively.}\label{fig:quadord}
\end{figure}
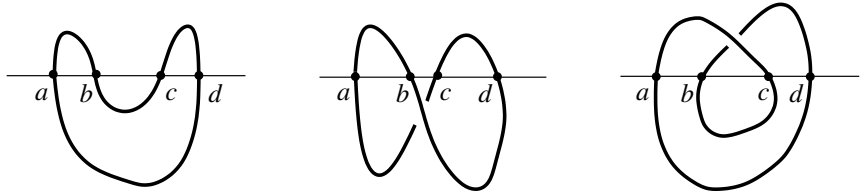

Similarly, the four points of a quadrisecant $abcd$ occur in that order
along the quadrisecant line, but may occur in any order along
the knot~$K$.  Of course, the order along~$K$ is only a cyclic order,
and if we ignore the orientation on~$K$ it is really just a dihedral order,
meaning one of the three cosets of the dihedral group~$D_4$ in~$S_4$.
Picking the lexicographically least element in each, we could
label these cosets $abcd$, $abdc$ and~$acbd$.  We will call
the corresponding classes of quadrisecants \emph{simple},
\emph{flipped} and \emph{alternating}, respectively, as in \figr{quadord}.
Note that this definition ignores the orientation of~$K$, and switching
the orientation of the quadrisecant does not change its type.

When discussing a quadrisecant $abcd$, we will usually
orient~$K$ so that~$b\in\arc{ad}$.  That means the cyclic order of
points along~$K$
will be $abcd$, $abdc$ or $acbd$, depending on the type of the quadrisecant.

In some sense, alternating quadrisecants are the most interesting.
These have also been called NSNS quadrisecants, because if
we view $\ray{cd}$ and $\ray{ba}$ as the North and South ends of
the midsegment $\segment{bc}$, then when $abcd$ is an alternating quadrisecant,
$K$ visits these ends alternately NSNS as it goes through the points
$acbd$.  It was noted by Cantarella, Kuperberg, Kusner and Sullivan
\cite{ckks} that the midsegment of an alternating
quadrisecant for $K$ is automatically in the \emph{second hull} of $K$.
Denne shows~\cite{denne} that nontrivial knots have alternating quadrisecants.
Budney, Conant, Scannell and Sinha~\cite{bcss} have shown that
the finite-type (Vassiliev) knot invariant of type $2$ can be computed
by counting alternating quadrisecants with appropriate multiplicity.

\section{Knots with unit thickness}
Because the ropelength problem is scale invariant, we find it
most convenient to rescale any knot $K$ to have thickness
(at least)~$1$.  This implies that $K$ is a $C^{1,1}$ curve
with curvature bounded above by~$2$. 

For any point $a\in \R^3$, let $B(a)$ denote the
open unit ball centered at~$a$.
Our first lemma, about the local structure of a thick knot, is by now standard.
(Compare \hbox{\cite[Lemma~4]{diao}} and \cite[Lemma~5]{cks}.)

\begin{lemma}\label{lem:thick}
Let $K$ be a knot of unit thickness.
If $a\in K$, then $B(a)$ contains a single unknotted arc
of~$K$; this arc has length at most~$\pi$ and is transverse to the
nested spheres centered at~$a$.
If $ab$ is a secant of~$K$ with $|a-b|<1$, then the ball of
diameter~$\segment{ab}$ intersects $K$ in a single unknotted arc
(either $\arc{ab}$
or $\arc{ba}$) whose length is at most $\arcsin|a-b|$.
\end{lemma}
\begin{proof}
If there were an arc of~$K$ tangent at some point $c$
to one of the spheres around~$a$ within~$B(a)$, then
triples near $(a,c,c)$ would have radius less than~$\half$.
For the second statement, if $K$ had a third intersection
point $c$ with the sphere of diameter $\segment{ab}$,
then we would have $r(a,b,c)<\half$.
The length bounds come from Schur's lemma.
\end{proof}

An immediate corollary is:
\begin{corollary}\label{cor:outside}
If $K$ has unit thickness, $a,b\in K$ and $p\in\arc{ab}$ with $a,b\notin B(p)$
then the complementary arc $\arc{ba}$ lies outside $B(p)$.
\end{corollary}

The following lemma should be compared to \cite[Lemma~5]{diao}
and \cite[Lemma~4]{cks}, but here we give a slightly
stronger version with a new proof.

\begin{lemma}\label{lem:balldist}
Let $K$ be a knot of unit thickness. If $a\in K$, then the radial
projection of~$K\setm\{a\}$ to the unit sphere $\bdry B(a)$ is $1$--Lipschitz,
that is, it does not increase length.
\end{lemma}

\begin{proof}
Consider what this projection does infinitesimally near a point $b\in K$.
Let $d=|a-b|$ and let $\theta$ be the angle at~$b$ between $K$ and the
segment $\segment{ab}$.  The projection stretches by a factor $1/d$ near $b$,
but does not see the radial part of the tangent vector to $K$.  Thus the
local Lipschitz constant on~$K$ is $(\sin\theta)/d$.
Now consider the circle through~$a$ and tangent to~$K$ at~$b$.
Plane geometry (see \figr{project}) shows
that its radius is $r=d/(2\sin\theta)$,
but it is a limit of circles through three points of~$K$, so
by the three-point characterization of thickness, $r$ is at least~$\tfrac12$;
that is, the Lipschitz constant~$1/2r$ is at most~$1$.
\end{proof}
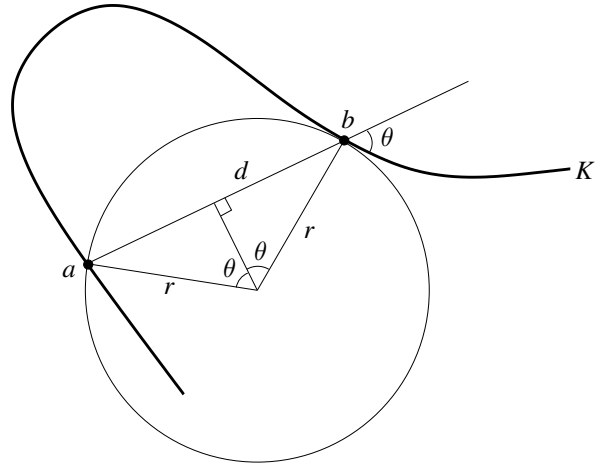
\begin{figure}[ht!]\centering
\begin{overpic}[scale=.4]{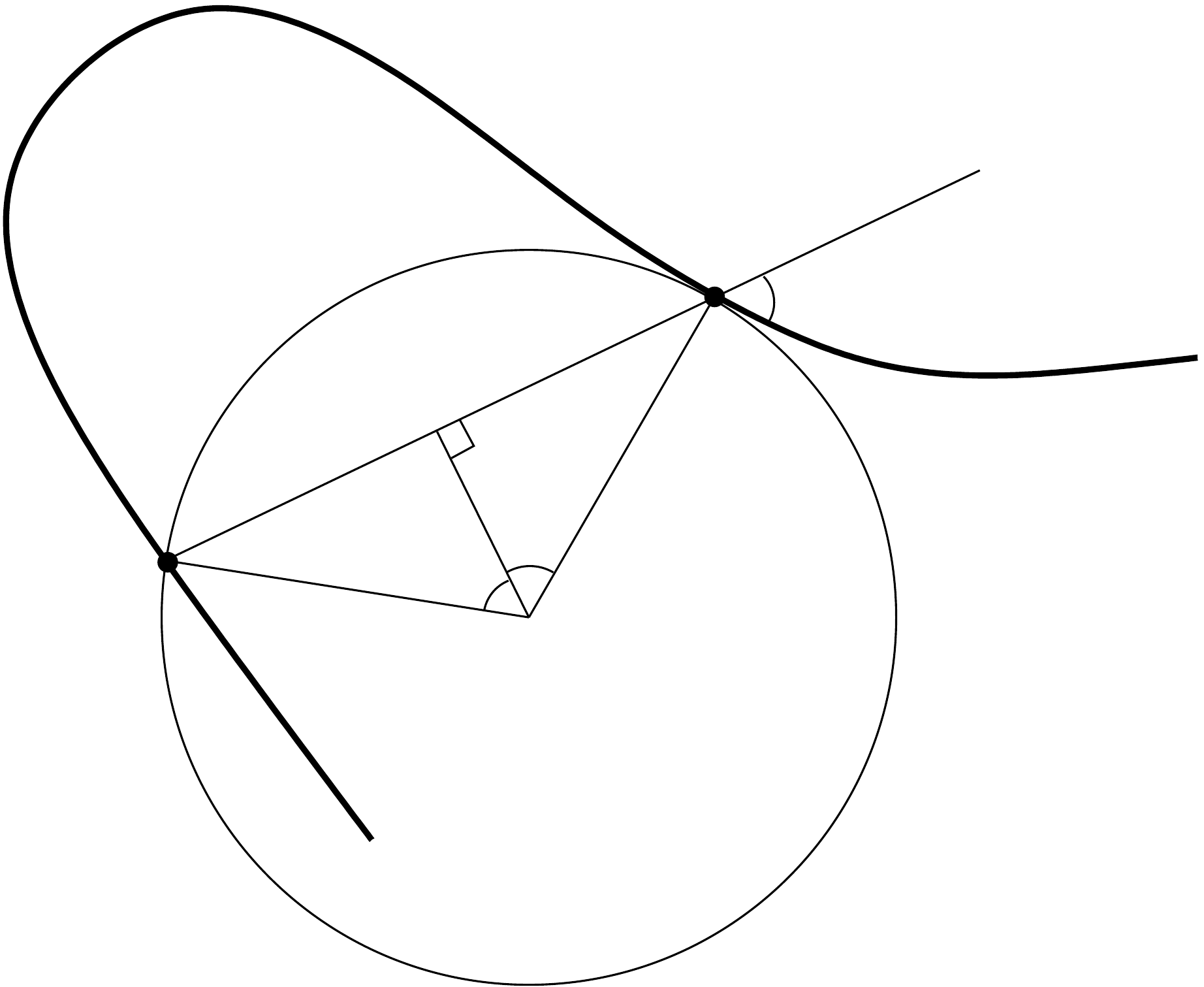}
{\small
\put(9,33){$a$}
\put(27,30){$r$}
\put(37.5,33){$\theta$}
\put(43,36){$\theta$}
\put(39.5,50.5){$d$}
\put(51.8,40){$r$}
\put(58.5,59.5){$b$}
\put(65.5,56.3){$\theta$}
\put(100,50.5){$K$}
}
\end{overpic}
\caption[Radius of a three-point circle]
{In the proof of \lem{balldist}, the circle tangent to~$K$ at~$b$
and passing through $a\in K$ has radius $r=d/(2\sin\theta)$
where $d=|a-b|$ and $\theta$ is the angle at~$b$ between~$K$ and~$\ray{ab}$.}
\label{fig:project}
\end{figure}

\begin{corollary}\label{cor:diao}
Suppose $K$ has unit thickness, and $p,a,b\in K$ with $p\notin \arc{ab}$.
Let $\angle apb$ be the angle between the vectors $a-p$ and $b-p$.
Then $\len{ab}\ge \angle apb$.
In particular, if~$apb$ is a reversed trisecant in $K$, then $\len{ab}\ge \pi$.
\end{corollary}
\begin{proof}
By \fullref{lem:balldist}, $\len{ab}$ is at least the length of the projected
curve on $B(p)$, which is at least the spherical distance
$\angle apb$ between its endpoints.  For a trisecant~$apb$ we have
$\angle apb=\pi$, and $p\notin \arc{ab}$ exactly when the trisecant is reversed.
\end{proof}

Observe that a quadrisecant $abcd$ includes four trisecants: $abc$,
$abd$, $acd$ and $bcd$.  Simple, flipped and alternating quadrisecants
have different numbers of reversed trisecants.  We can apply \cor{diao}
to these trisecants to get simple lower bounds on ropelength for any
curve with a quadrisecant.

\begin{theorem}\label{thm:quadbd}
The ropelength of a knot with a simple, flipped or alternating quadrisecant
is at least $\pi$, $2\pi$ or $3\pi$, respectively.
\end{theorem}
\begin{proof}
Rescale $K$ to have unit thickness, so that its ropelength equals its
length~$\ell$.
Let~$abcd$ be the quadrisecant, and orient~$K$ in the usual way,
so that $b\in\arc{ad}$.
In the case of a simple quadrisecant,
the trisecant $dba$ is reversed,
so $\ell\ge\len{da}\ge \pi$, using \cor{diao}.
In the case of a flipped quadrisecant,
the trisecants $cba$ and $bcd$ are reversed,
so $\ell\ge\len{ca}+\len{bd}\ge2\pi$.
In the case of an alternating quadrisecant,
the trisecants $abc$, $bcd$ and $dca$ are reversed,
so $\ell\ge\len{ac}+\len{bd}+\len{da}\ge3\pi$.
\end{proof}

Of course, any closed curve has ropelength at least~$\pi$, independent
of whether it is knotted or has any quadrisecants, because its total curvature
is at least $2\pi$.  But curves arbitrarily
close (in the $C^1$ or even $C^\infty$ sense) to a round circle can have
simple quadrisecants, so at least the first bound in the theorem above is sharp.

Although Kuperberg has shown that any nontrivial (tame)
knot has a quadrisecant, and
Denne shows that in fact it has an alternating quadrisecant, the
bounds from \thm{quadbd} are not as good as the previously known
bounds of~\cite{diao} or even~\cite{cks}.  To improve our bounds, in \secn{ess}
we will consider Kuperberg's notion of an essential quadrisecant.

\section{Length bounds in terms of segment lengths}
Given a thick knot~$K$ with quadrisecant $abcd$, we can bound its
ropelength in terms of the distances along the quadrisecant line.
Whenever we discuss such a quadrisecant, we will abbreviate these
three distances as $r := |a-b|$, $s := |b-c|$ and $t := |c-d|$.
We start with some lower bounds for $r$, $s$ and~$t$ for
quadrisecants of certain types.

\begin{lemma}\label{lem:flipseg}
  If $abcd$ is a flipped quadrisecant for a knot of unit thickness,
  then the midsegment has length $s\ge 1$.  Furthermore if $r\ge 1$
  then the whole arc $\arc{ca}$ lies outside~$B(b)$; similarly if $t\ge 1$
  then $\arc{bd}$ lies outside $B(c)$.
\end{lemma}
\begin{proof}
Orient the knot in the usual way.  If $s=|b-c|<1$, then by \lem{thick}
either $\len{cab}<\pi/2$ or $\len{bdc}<\pi/2$.
But since $cba$ and $bcd$ are reversed trisecants, we have
$\len{ca}\ge\pi$ and $\len{bd}\ge\pi$.  This is a contradiction
because $\len{cab}=\len{ca}+\len{ab}$ and~$\len{bdc}=\len{bd}+\len{dc}$.
The second statement follows immediately from \cor{outside}.
\end{proof}

\begin{lemma}\label{lem:altarc}
If $abcd$ is an alternating quadrisecant for a knot of unit thickness,
then $r\ge1$ and $t\ge1$.
With the usual orientation, the entire arc~$\arc{da}$
thus lies outside $B(b)\cup B(c)$.
If $s\ge1$ as well, then $\arc{ac}$ lies outside~$B(b)$
and $\arc{bd}$ lies outside $B(c)$.
\end{lemma}
\begin{proof}
If $r=|a-b|<1$, then by \lem{thick} either $\len{acb}<\pi/2$
or $\len{bda}<\pi/2$.
Similarly, if $t=|c-d|<1$, then either $\len{cbd}<\pi/2$ or $\len{dac}<\pi/2$.
But as in the proof of \thm{quadbd}, we have
$\len{ac}\ge \pi$ and $\len{bd}\ge \pi$, contradicting any choice of
the inequalities above.  Thus we have $r,t\ge1$.
Because $a$ and $d$ are outside $B(b)$ and~$B(c)$, the remaining
statements follow from \cor{outside}.
\end{proof}

As suggested by the discussion above,
we will often find ourselves in the situation where
we have an arc of a knot known to stay outside a unit ball.
We can compute exactly the minimum length of such an arc
in terms of the following functions.

\begin{definition}
For $r\ge 1$, let $f(r) := \sqrt{r^2-1}+\arcsin(1/r)$.
For $r,s\ge 1$ and $\theta\in[0,\pi]$, the minimum length function is defined by
$$ m(r,s,\theta) := \left\{ \begin{array}{ll}
  \sqrt{r^2+s^2-2rs\cos\theta}
      &\quad\textrm{if } \theta \le \arccos(1/r)+\arccos(1/s), \\
  f(r)+f(s)+(\theta-\pi)
      &\quad\textrm{if } \theta \ge \arccos(1/r)+\arccos(1/s).
\end{array}\right.$$
\end{definition}

The function $f(r)$ will arise again in other situations.
The function $m$ was defined exactly to make the following bound sharp:

\begin{lemma}\label{lem:min}
Any arc $\gamma$ from $a$ to $b$, staying outside $B(p)$,
has length at least ${m(|a-p|,|b-p|,\angle apb)}$.
\end{lemma}
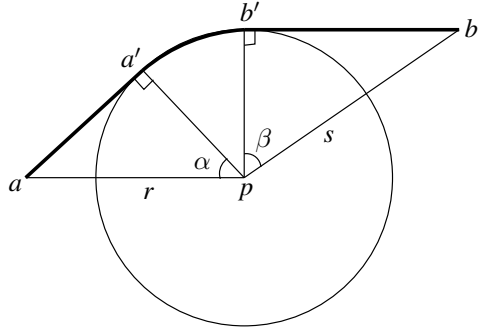
\begin{figure}[ht!]\centering
\begin{overpic}[scale=.5]{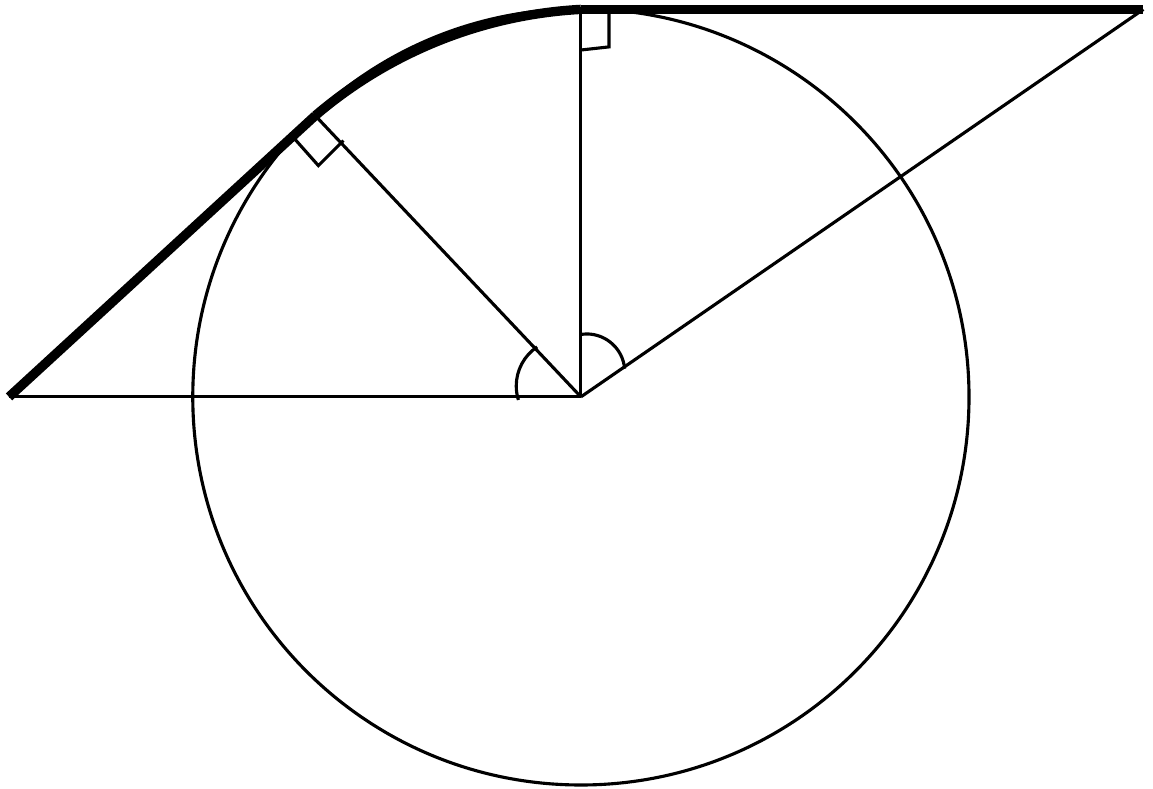}
{\small
\put(-3.5,31){$a$}
\put(27,29){$r$}
\put(38.5,36){$\alpha$}
\put(52.5,41){$\beta$}
\put(49,30){$p$}
\put(22,59){$a'$}
\put(68,42){$s$}
\put(49,69.5){$b'$}
\put(100,66){$b$}
}
\end{overpic}
\caption[Shortest curve avoiding a ball]
{If points $a,b$ are at distances $r,s\ge 1$ (respectively) from~$p$,
then the shortest curve from~$a$ to~$b$ avoiding $B(p)$ is planar.  Either
it is a straight segment or (in the case illustrated) it includes an arc
of $\bdry B(p)$.  In either case, its length is $m(r,s,\angle{apb})$.}
\label{fig:minlen}
\end{figure}

\begin{proof}
Let $r:=|a-p|$ and $s:=|b-p|$ be the distances to $p$ (with $r,s\ge 1$)
and let $\theta:=\angle apb$ be the angle between $a-p$ and $b-p$.
The shortest path from $a$ to $b$ staying outside $B:=B(p)$
either is the straight segment or is the $C^1$ join of a straight
segment from $a$ to $\bdry B$, 
a great-circle arc in $\bdry B$, and a straight segment from $\bdry B$ to~$b$.
In either case, we see that the path lies in the plane through $a$, $p$, $b$
(shown in \figr{minlen}).
In this plane, draw the lines from~$a$ and~$b$
tangent to~$\bdry B$.
Let $\alpha:=\angle apa'$ and $\beta:=\angle bpb'$, where~$a'$ and~$b'$ are
the points of tangency.  Then $\cos\alpha=1/r$ and $\cos\beta=1/s$.
Clearly if $\alpha +\beta\ge\theta$ then the shortest path is the
straight segment from~$a$ to~$b$, with length $\sqrt{r^2+s^2-2rs\cos\theta}$.
If $\alpha +\beta\le\theta$ then the 
shortest path consists of the $C^1$ join described above, with length
$$\sqrt{r^2-1} + \big(\theta - (\alpha +\beta)\big) + \sqrt{s^2-1}
= f(r)+f(s)+(\theta-\pi).\proved$$
\end{proof}

An important special case is when $\theta=\pi$.  Here we are
always in the case $\alpha+\beta\le\theta$, so we get the following
corollary.
\begin{corollary}\label{cor:mindist}
If $a$ and~$b$ lie at distances $r$ and~$s$ along opposite rays
from $p$ (so that~$\angle apb = \pi$) then the length of any arc
from $a$ to $b$ avoiding $B(p)$ is at least
$$f(r)+f(s) =
    \sqrt{r^2-1} + \arcsin(1/r) + \sqrt{s^2-1} + \arcsin(1/s).$$
\end{corollary}

We note that the special case of this formula when $r=s$ also appears
in recent papers by Dumitrescu, Ebbers-Baumann, Gr\"une, Klein and Rote
\cite{degkr,egk} investigating the geometric dilation (or distortion)
of planar graphs.  Gromov had given a lower bound for the distortion of
a closed curve (see the paper by Kusner and Sullivan \cite{KS-disto});
in~\hbox{\cite{degkr,egk}} sharper bounds in terms of the diameter and width
of the curve are derived using this minimum-length arc avoiding a ball.
(Although the bounds are stated there only for plane curves they apply
equally well to space curves.)  See also the article by Denne and Sullivan
\cite{DS-disto} for further development of these ideas and a proof that
knotted curves have distortion more than~$4$.

\begin{lemma}\label{lem:longarc}
Let $abcd$ be an alternating quadrisecant for a knot of unit
thickness (oriented in the usual way).
Let $r:=|a-b|$, $s:=|b-c|$ and $t:=|c-d|$ be the lengths of the segments
along $abcd$.  Then $\len{ad} \ge f(r) + s + f(t)$.
The same holds if $abcd$ is a simple quadrisecant as long as $r,t\ge1$.
\end{lemma}
\begin{proof}
  In either case (and as we already noted in \lem{altarc}
  for the alternating case) we find, using \cor{outside}, that
  $\arc{da}$ lies outside $B(b)\cup B(c)$.
  As in the proof of \lem{min}, the shortest arc from~$d$ to~$a$
  outside these balls will be the $C^1$ join of various pieces:
  these alternate between straight segments in space and great-circle
  arcs in the boundaries of the balls.
  Here, the straight segment in the middle always has length exactly
  $s:=|b-c|$. As in \cor{mindist}, the overall length is then at
  least~$f(r)+s+f(t)$ as desired.
\end{proof}

\section{Essential secants}\label{sec:ess}
We have seen that the existence of a quadrisecant for~$K$ is not enough to
get good lower bounds on ropelength, because some quadrisecants do
not capture the knottedness of~$K$.  Kuperberg~\cite{kup} introduced
the notion of \emph{essential} secants and quadrisecants
(which he called ``topologically nontrivial'').  We will see below
that these give us much better ropelength bounds.

We extend Kuperberg's definition to say when an arc $\arc{ab}$ of a knot~$K$
is essential, capturing part of the knottedness of $K$.
Generically, the knot~$K$ together with the segment~$S=\segment{ab}$
forms a knotted $\Theta$--graph in space (that is, a graph with three
edges connecting the same two vertices).  To adapt 
Kuperberg's definition, we consider such knotted $\Theta$--graphs.

\begin{figure}[ht!]\centering
\begin{overpic}[scale=.4]{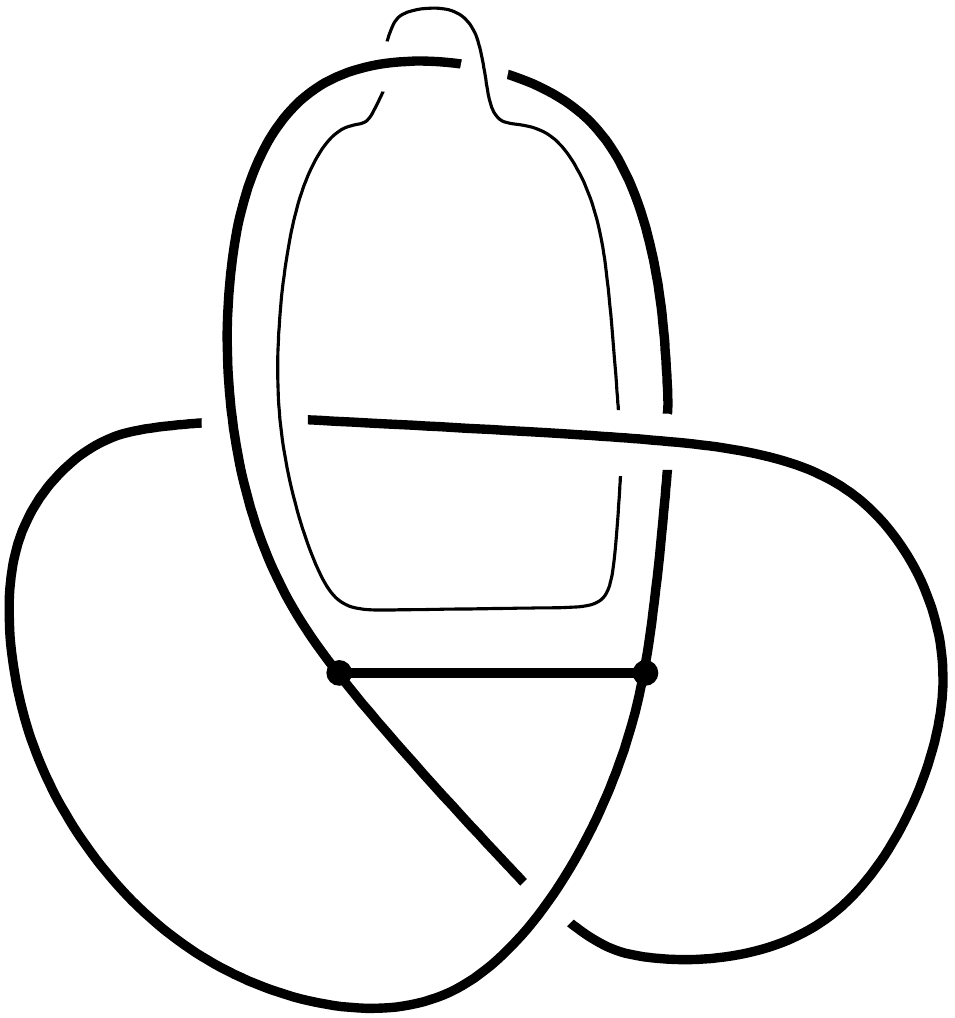}
{\small
\put(21.5,90){$\alpha$}
\put(24.7,31.2){$p$}
\put(66,31.5){$q$}
\put(47,26){$\beta$}
\put(38,6){$\gamma$}
\put(52,74){$\delta$}
}
\end{overpic}
\caption[An essential arc]
{In this knotted $\Theta$--graph $\alpha\cup\beta\cup\gamma$, the
ordered triple $(\alpha,\beta,\gamma)$ is essential.  To see this, we find
the parallel~$\delta=h(\alpha,\beta,\gamma)$ to $\alpha\cup\beta$
which has linking number zero with $\alpha\cup\gamma$, and note that it is
homotopically nontrivial in the knot complement $\R^3\setm(\alpha\cup\gamma)$.
In this illustration, $\beta$ is the straight segment~$\segment{pq}$, so
we equally say that the arc~$\alpha$ of the knot $\alpha\cup\gamma$
is essential.}
\label{fig:essdef} \end{figure}
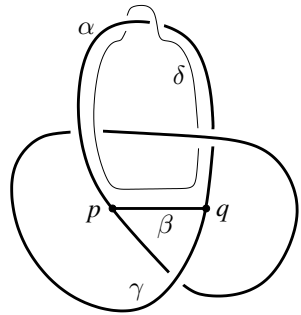

\begin{definition}
  Suppose $\alpha$, $\beta$ and $\gamma$ are three disjoint simple
  arcs from~$p$ to~$q$, forming a knotted $\Theta$--graph.
  Then we say that the ordered triple $(\alpha,\beta,\gamma)$
  is \emph{inessential}
  if there is a disk~$D$ bounded by the knot $\alpha\cup\beta$ having no
  interior intersections with the knot~$\alpha\cup\gamma$.
  (We allow self-intersections of~$D$, and interior intersections with $\beta$;
  the latter are certainly necessary if $\alpha\cup\beta$ is knotted.)

  An equivalent definition, illustrated in \figr{essdef}, is as follows:
  Let $X:=\R^3\setm(\alpha\cup\gamma\big)$,
  and consider a parallel curve~$\delta$ to $\alpha\cup\beta$ in~$X$.
  Here by \emph{parallel} we mean that $\alpha\cup\beta$
  and~$\delta$ cobound an annulus embedded in~$X$.
  We choose the parallel to be homologically trivial in~$X$.
  (Since the homology of the knot complement~$X$ is~$\Z$,
  this simply means we take~$\delta$ to have linking number
  zero with $\alpha\cup\gamma$.  This determines $\delta$ uniquely
  up to homotopy.)
  Let $x_0\in\delta$ near $p$ be a basepoint for~$X$,
  and let $h=h(\alpha,\beta,\gamma)\in\pi_1(X,x_0)$
  be the homotopy class of~$\delta$.
  Then $(\alpha,\beta,\gamma)$ is \emph{inessential} if $h$ is trivial.

  We say that $(\alpha,\beta,\gamma)$ is \emph{essential}
  if it is not inessential, meaning that $h(\alpha,\beta,\gamma)$ is nontrivial.

  Now let $\lambda$ be a meridian loop (linking $\alpha\cup\gamma$ near~$x_0$)
  in the knot complement~$X$.
  If the commutator $[h(\alpha,\beta,\gamma),\lambda]$ is nontrivial
  then we say $(\alpha,\beta,\gamma)$ is \emph{strongly essential}.
\end{definition}

This notion is clearly a topological invariant of the (ambient isotopy)
class of the knotted $\Theta$--graph.
For an introduction to the theory of knotted graphs, see the article by
Kauffman \cite{kauf89};
note that since the vertices of the $\Theta$--graph have degree three, in
our situation there is no distinction between what Kauffman calls
topological and rigid vertices.
The three arcs approaching one vertex can be braided arbitrarily without
affecting the topological type of the knotted $\Theta$--graph.

\begin{lemma}\label{lem:stress}
In a knotted $\Theta$--graph $\alpha\cup\beta\cup\gamma$, the
triple $(\alpha,\beta,\gamma)$ is strongly essential if and only if
$(\gamma,\beta,\alpha)$ is.
\end{lemma}
\begin{proof}
The homotopy classes $h=h(\alpha,\beta,\gamma)$ and $h'=h(\gamma,\beta,\alpha)$
proceed outwards from~$x_0$ along~$\beta$ and then
return backwards along~$\alpha$ or~$\gamma$.
Note that the product $h^{-1} h'$
is homotopic to a parallel of the knot $\alpha\cup\gamma$.
Since a torus has abelian fundamental group,
this parallel commutes with the meridian~$\lambda$.
It follows that $[\lambda,h]=[\lambda,h']$.
\end{proof}

\begin{figure}[ht!]\centering
\begin{overpic}[scale=.5]{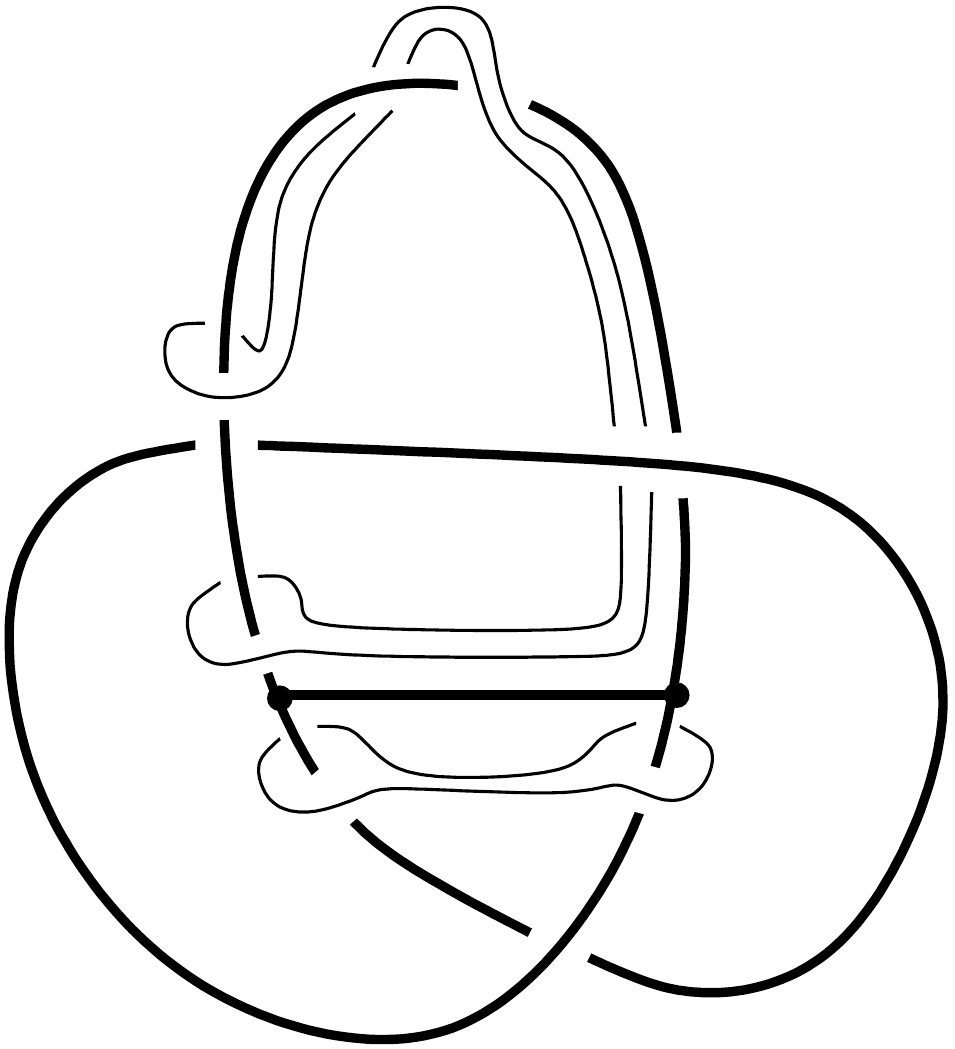}
{\small
\put(20.3,31.2){$p$}
\put(67,31.5){$q$}
\put(22,90){$\alpha$}
\put(50,94){$[\lambda,h(\alpha,\beta)]$}
\put(44,26.2){${}^\beta$}
\put(21.7,19.2){$l_{\beta}$}
\put(8.7,5.7){$\gamma$}
}
\end{overpic}
\caption[Loop around a secant segment]
{If $\lambda$ is a meridian curve linking $\alpha\cup\gamma$, then
the commutator $[\lambda,h(\alpha,\beta)]$ is homotopic
to the loop~$l_{\beta}$ along~$\beta$.} 
\label{fig:strongess}
\end{figure}
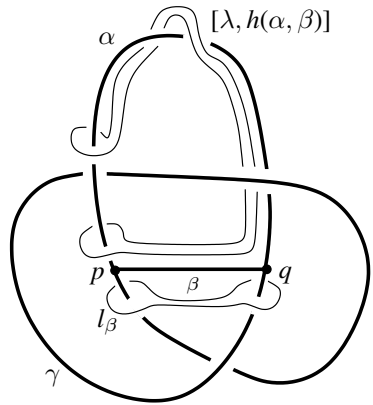

This commutator $[\lambda,h(\alpha,\beta,\gamma)]$ that comes up in the
definition of strongly essential will later be referred to as the
\emph{loop~$l_\beta$ along~$\beta$}; it can be represented by a curve which
follows a parallel~$\beta'$ of~$\beta$, then loops around $\alpha\cup\gamma$
along a meridian near~$q$, then follows~$\beta'^{-1}$, then loops backwards
along a meridian near~$p$, as in \figr{strongess}.

We apply the definition of essential to arcs of a knot as follows.

\begin{definition}
  If $K$ is a knot and $a,b\in K$, let $S=\segment{ab}$.
  We say $\arc{ab}$ is \emph{(strongly) essential} in~$K$
  if for every $\eps>0$ there exists some $\eps$--perturbation of~$S$
  (with endpoints fixed) to a curve $S'$ such that
  $K\cup S'$ is an embedded $\Theta$--graph in which
  $(\arc{ab},S',\arc{ba})$ is (strongly) essential.
\end{definition}

\begin{remark}
Allowing the $\eps$--perturbation
ensures that the set of essential secants is closed in the set of all secants of~$K$,
and lets us handle the case when $S$ intersects~$K$.
We could allow the perturbation only in that case of intersection;
the combing arguments of~\cite{DS-conv} show
the resulting definition is equivalent.
We require only that $S'$ be~$\eps$--close to~$S$ in the $C^0$ sense; it thus
could be locally knotted, but in the end we care only about
the homotopy class~$h$, and not an isotopy class.
\end{remark}

In~\cite{ckks} it was shown that if $K$ is an unknot, then any arc
$\arc{ab}$ is inessential.  In our context, this follows immediately,
because the homology and homotopy groups of $X:=\R^3\setm K$ are equal for
an unknot, so any curve $\delta$ having linking number zero
with $K$ is homotopically trivial in $X$.
We can use Dehn's lemma to prove a converse statement:
\begin{theorem}
If $a,b\in K$ and both $\arc{ab}$ and $\arc{ba}$ are inessential,
then $K$ is unknotted.
\end{theorem}
\begin{proof}
Let $S$ be the secant $\segment{ab}$, perturbed if necessary to avoid
interior intersections with~$K$.
We know that $\arc{ab}\cup S$ and $\arc{ba}\cup S$
bound disks whose interiors are disjoint from~$K$.
Glue these two disks together along~$S$ to form a disk~$D$ spanning~$K$.
This disk may have self intersections, but these occur away
from~$K$, which is the boundary of~$D$.
By Dehn's lemma, we can replace~$D$ by an embedded disk, hence $K$ is unknotted.
\end{proof}

\begin{definition}
  A secant $ab$ of $K$ is \emph{essential} if
  both subarcs $\arc{ab}$ and $\arc{ba}$ are essential.
  A secant $ab$ is \emph{strongly essential} if $\arc{ab}$
  (or, equivalently, $\arc{ba}$) is strongly essential.
\end{definition}

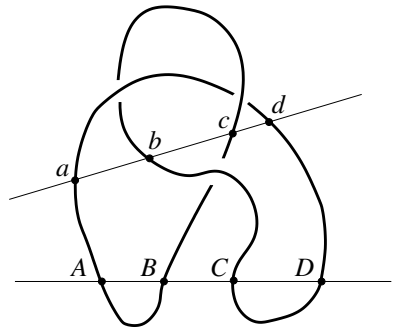
\begin{figure}[ht!]\centering
\begin{overpic}[scale=.4]{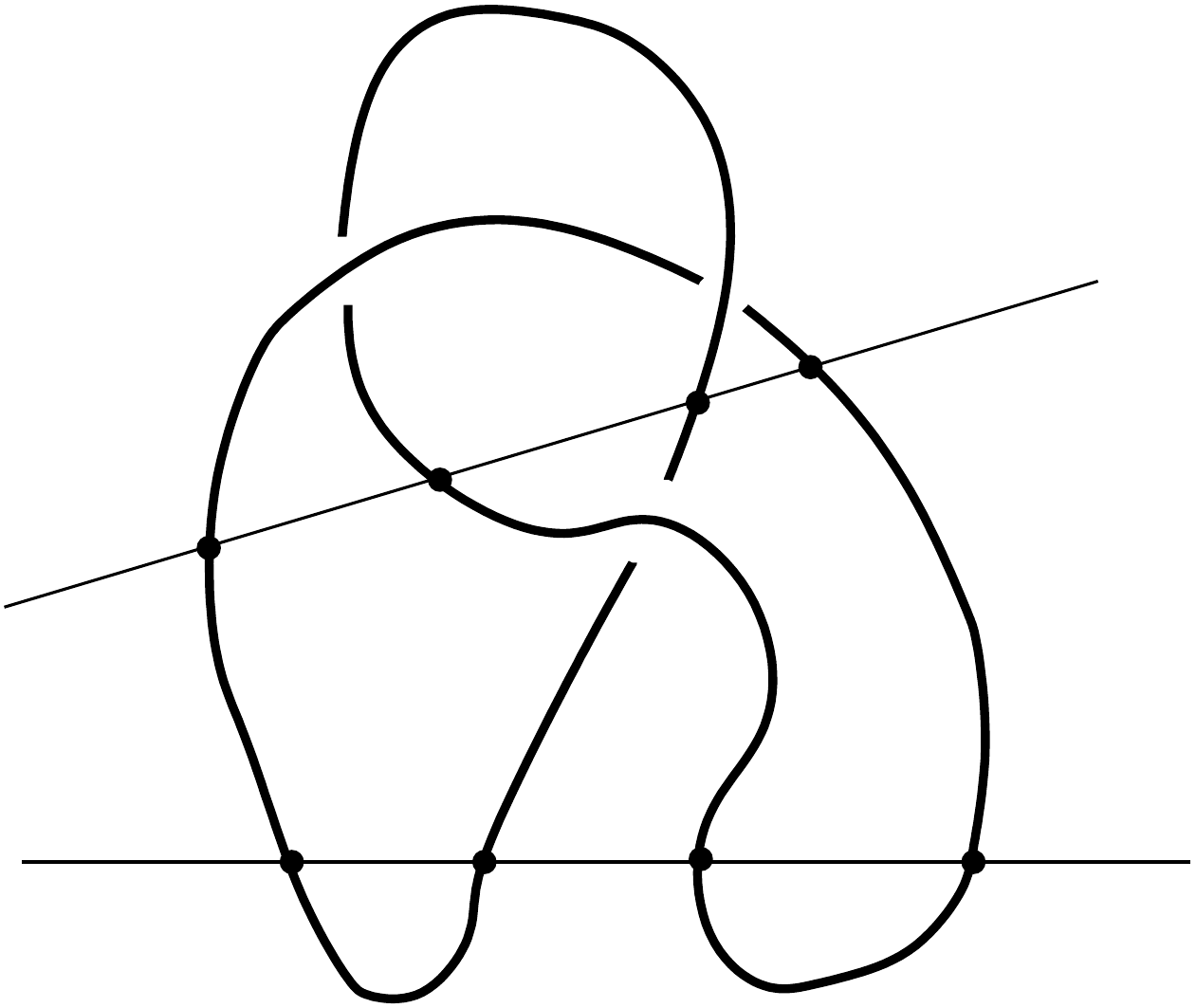}
{\small
\put(12,39.5){$a$}
\put(36,46.5){$b$}
\put(54,52){$c$}
\put(68,56){$d$}
\put(16,13.5){$A$}
\put(34,13.5){$B$}
\put(52,13.5){$C$}
\put(74,13.5){$D$}
}
\end{overpic}
\caption[Essential and inessential quadrisecants]
{This trefoil knot has two quadrisecants. Quadrisecant $abcd$ is
alternating and essential (meaning that $bc$ is essential, although
here in fact also $ab$ and $cd$ are essential).
Quadrisecant $ABCD$ is simple and inessential,
since $AB$ and $CD$ are inessential
(although $BC$ is essential).}\label{fig:twoquad}
\end{figure}

To call a quadrisecant $abcd$ essential, we could follow Kuperberg and
require that the secants $ab$, $bc$ and $cd$ all be essential.
But instead, depending on the order type of the quadrisecant,
we require this only of those secants whose length could not already
be bounded as in \thm{quadbd}, namely those secants whose endpoints
are consecutive along the knot.
That is, for simple quadrisecants, all three secants must be essential;
for flipped quadrisecants the end secants $ab$ and $cd$ must be essential;
for alternating quadrisecants, the middle secant $bc$ must be essential.
\figr{twoquad} shows a knot with essential and inessential quadrisecants.

More formally, we can give the following definition for any $n$--secant.
\begin{definition}
An $n$--secant $a_1a_2\ldots a_n$ is \emph{essential} if
we have $a_ia_{i+1}$ essential for each~$i$ such that one of the arcs
$\arc{a_ia_{i+1}}$ and $\arc{a_{i+1}a_i}$ includes no other~$a_j$.
\end{definition}

Kuperberg introduced the notion of essential secants and showed the
following result.
\begin{theorem}\label{thm:kup}
If $K$ is a nontrivial knot parameterized by a generic polynomial,
then~$K$ has an essential quadrisecant.
\end{theorem}

As mentioned above, Kuperberg did not distinguish the different order types,
and so he actually obtained a quadrisecant $abcd$ where all three segments $ab$,
$bc$ and $cd$ are essential.  Since this is more than we will need here,
we have given our weaker definition of essential for
the flipped and alternating cases, making these easier to produce.

Kuperberg used the fact that a limit of essential quadrisecants must
still be a quadrisecant in order to show that every knot has a
quadrisecant.  Since we want an essential quadrisecant for every knot,
we next need to show that being essential is preserved in such limits.

\section{Limits of essential secants}
Being essential is a topological property of a knotted $\Theta$--graph.
One approach to show that a limit of essential secants remains essential
is to show that nearby knotted $\Theta$s are isotopic.  In fact,
we have pursued this approach in another paper~\cite{DS-conv}
where we show that given any knotted graph of finite total curvature,
any other graph which is sufficiently close (in a $C^1$ sense) is isotopic.
(Our definition makes essential a closed condition, so the case where
a secant has interior intersections with the knot,
and is thus not a knotted theta, causes no trouble.)

Here, however, our knots are thick, hence $C^{1,1}$, and we give a
simpler direct argument for the limit of essential secants.

\begin{lemma}\label{lem:ambiso}
Let $K$ be a knot of thicknes~$\tau>0$, and $K'$ be a $C^1$ knot which is
close to~$K$ in the following sense: corresponding points $p$ and~$p'$
are within distance $\eps<\tau/4$ and their tangent vectors are within angle $\pi/6$. 
Then $K$ and~$K'$ are ambient isotopic; the isotopy can be chosen to
move each point by a distance less than~$\eps$.
\end{lemma}
Note that the constants here are sharp within a factor of two or three:
if the distance from $K$ to~$K'$ exceeds $\tau/2$ we can have
strand passage, while if the angle between tangent vectors
exceeds $\pi/2$ we can have local knotting in $K'$.
\begin{proof}
Rescale so $K$ has unit thickness and $\eps<1/4$.
Clearly $K'$ lies within the thick tube around~$K$.
Each point $p'\in K'$ corresponds to some point $p\in K$,
but also has a unique nearest point $p_0\in K$, which is within distance
$|p'-p|<\eps$ of~$p'$, hence within~$2\eps$ of $p$.
By \lem{thick}, the arclength from $p_0$
to~$p$ is at most $\arcsin2\eps$, so the angle between the tangent
vectors there is at most $2\arcsin2\eps<\pi/3$.
The point~$p'$ is thus in the normal disk at~$p_0$ and has tangent
vector within $\pi/2$ of that at~$p_0$.  In other words, $K'$
is transverse to the foliation of the thick tube by normal disks.
Construct the isotopy from~$K'$ to~$K$ as the union of isotopies in
these disks: On each disk, we move $p'$ to~$p_0$, coning this outwards
to the fixed boundary.  No other point moves further than~$p'$
which moves at most distance~$\eps$.
\end{proof}

\begin{proposition}\label{pr:limess}
If the $C^{1,1}$ knots $K_i$ have essential arcs $\arc{a_ib_i}$,
and if the~$K_i$ converge in~$C^1$ to some thick limit knot~$K$,
with $a_i\to a$ and $b_i\to b$, then the arc $\arc{ab}$ is essential for~$K$.
\end{proposition}

\begin{proof}
We can reduce to the case $a_i=a$, $b_i=b$ by
applying euclidean similarities (approaching the identity) to the $K_i$.

Given any $\eps>0$, we prove there is a $2\eps$--perturbation of $\segment{ab}$
making an essential knotted $\Theta$.  Then by definition, $\arc{ab}$ is essential.

For large enough~$i$, the knot $K_i$ is within $\eps$ of $K$.
Let~$I_i$ be the ambient isotopy described in \lem{ambiso} with $K=I_i(K_i)$.
Since $K_i$ is essential, by definition,
we can find an $\eps$--perturbation $S'_i$ of~$S_i$
such that $\Theta_i := K_i\cup S'_i$ is an embedded essential $\Theta$--graph.
Setting $S''_i:=I_i(S'_i)$, this is the
desired $2\eps$ perturbation of $\segment{ab}$.
By definition $K_i\cup S'_i$ is isotopic to $K\cup S''_i$,
so the latter is also essential.
\end{proof}

\begin{corollary}\label{cor:alless}
Every nontrivial $C^{1,1}$ knot has an essential quadrisecant.
\end{corollary}
\begin{proof}
Any $C^{1,1}$ knot~$K$ is a $C^1$--limit of polynomial knots,
which can be taken to be generic in the sense of~\cite{kup}.
Then by \thm{kup} they have essential quadrisecants.
Some subsequence of these quadrisecants converges to a
quadrisecant for~$K$, which is essential by \prop{limess}.
\end{proof}

\section{Arcs becoming essential}
We showed in \cor{alless} that every nontrivial $C^{1,1}$ knot has
an essential quadrisecant. Our aim is to find the least length of an
essential arc $\arc{pq}$ for a thick knot~$K$, and use this to
get better lower bounds on the ropelength of knots.
This leads us to consider what happens when arcs change from inessential to essential.

\begin{theorem}\label{thm:change} 
Suppose $\arc{ac}$ is in the boundary of the set of essential arcs
for a knot~$K$.  (That is, $\arc{ac}$ is essential, but there
are inessential arcs of $K$ with endpoints arbitrarily close to~$a$ and~$c$.)
Then $K$ must intersect the interior of segment~$\segment{ac}$,
and in fact there is some essential trisecant $abc$.
\end{theorem}

\begin{proof}
If $K$ did not intersect the interior of segment~$\segment{ac}$
(in a component separate from~$a$ and~$c$)
then all nearby secant segments would form ambient isotopic $\Theta$--graphs,
and thus would all be essential.
Therefore $a$ and~$c$ are the first and third points of some trisecant~$abc$.
Let $S$ and~$S'$ be two perturbations of~$\segment{ac}$,
forming $\Theta$--graphs with~$K$, such that
$(\arc{ac},S,\arc{ca})$ is essential but $(\arc{ac},S',\arc{ca})$ is not.
(The first exists because~$\arc{ac}$ is essential, and the second because it
is near inessential arcs.)

For clarity, we first assume that $b$ is the only point
of~$K$ in the interior of $\segment{ac}$ and that
$S$ and~$S'$ differ merely by going to the
two different sides of~$K$ near~$b$, as in \figr{secperturb}.
(We will return to the more general case at the end of the proof.)
We will show that $ab$ is a strongly essential secant; by symmetry the same
is true of~$bc$, and thus $abc$ is an essential trisecant, as desired.

\begin{figure}[ht!]\centering
\begin{overpic}[scale=.5]{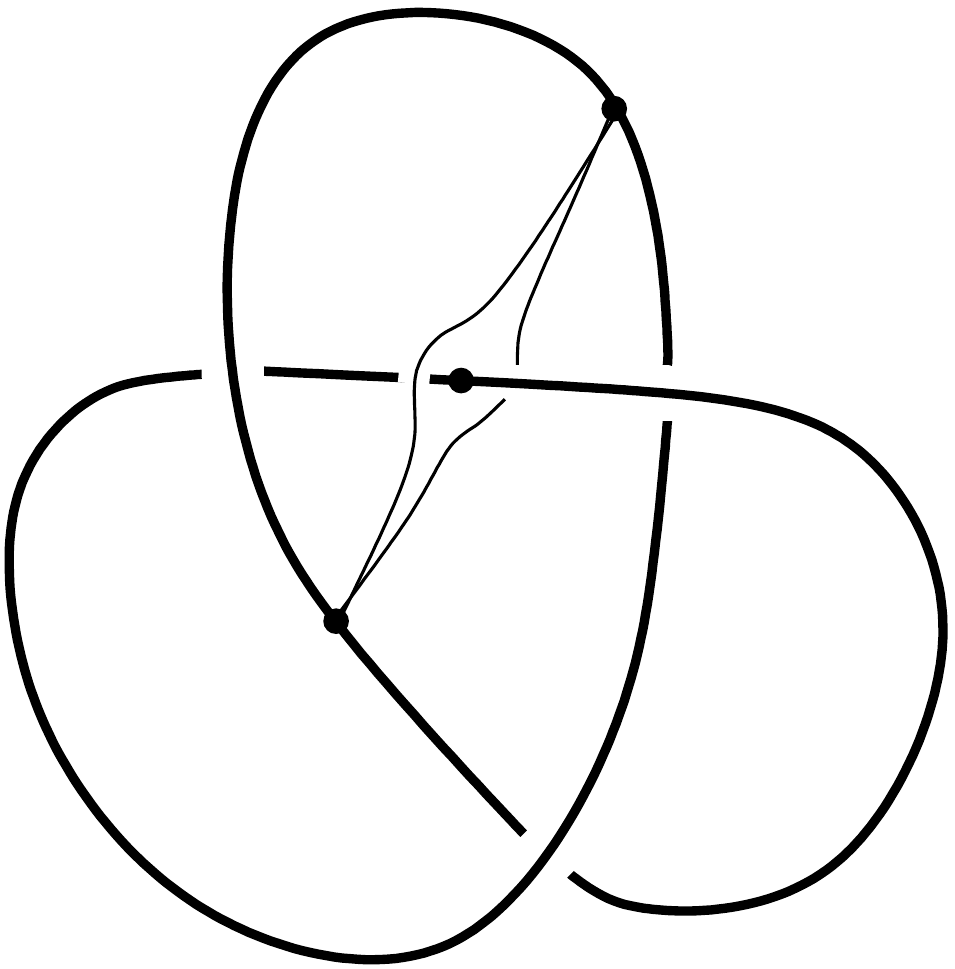}
{\small
\put(29,32){$a$}
\put(50,50){$b$}
\put(64,91){$c$}
\put(48,75){$S'$}
\put(58,72){$S$}
}
\end{overpic}
\caption[A borderline essential secant]
{Secant~$ac$ is essential but some nearby secants are not.
By \thm{change} there must be an essential trisecant $abc$,
because there are perturbations $S$ and $S'$ of $\segment{ac}$ which are
essential and inessential, respectively.}
\label{fig:secperturb}
\end{figure}
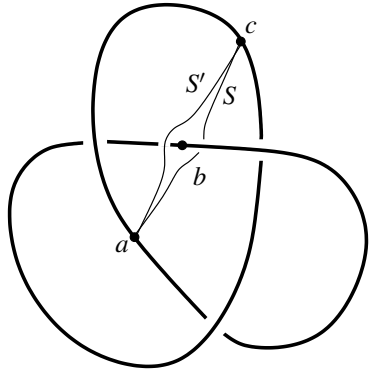

By the definition of essential, the homotopy class
$h:=h(\arc{ac},S,\arc{ca})$ is nontrivial in~$\pi_1(\R^3\setm K)$,
but $h':=h(\arc{ac},S',\arc{ca})$ is trivial.
Since both have linking number zero with~$K$,
they differ not only by the meridian loop around~$K$ near~$b$
(seen in the change from~$S$ to~$S'$) but also by a meridian
loop around~$K$ somewhere along the arc~$\arc{ac}$, say near~$a$.
Let $\delta$ and~$\delta'$ be the standard loops representing
these homotopy classes (as in the definition of essential), and
consider the subarc of $\delta$ which follows along~$\segment{bc}$ and
then back along~$\arc{ac}$.
The fact that $\delta'$ is null-homotopic means that this subarc
is homotopic to a parallel to $\segment{ba}$. This means $\delta$ is homotopic
to the loop $l_{\segment{ab}}$ along $\segment{ab}$, as in \figr{loop}.
Thus $l_{\segment{ab}}$ represents the nontrivial homotopy class $h$,
so by definition~$\segment{ab}$ is strongly essential.

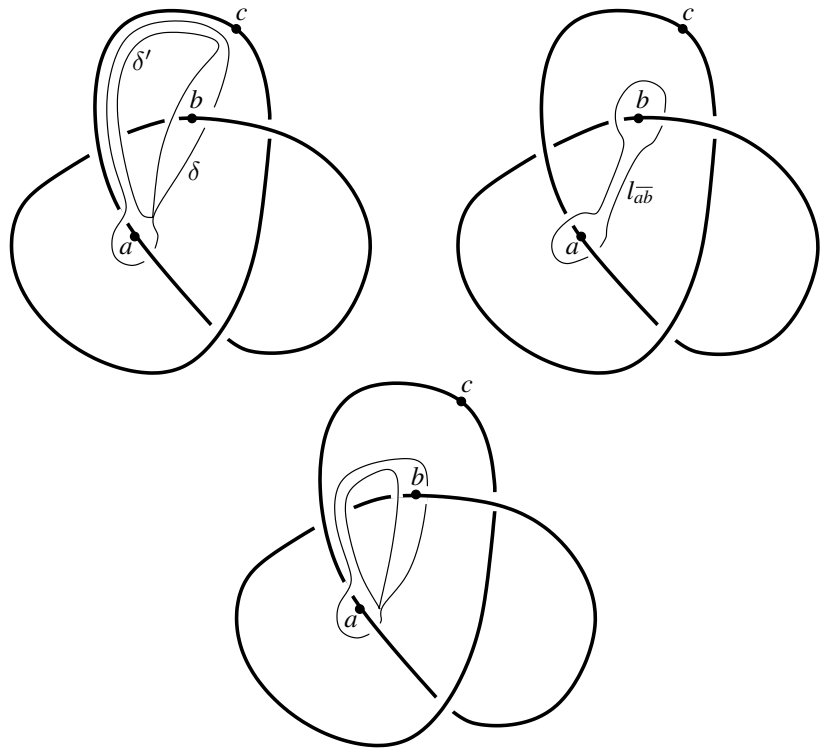
\begin{figure}[ht!]\small\centering
\centerline{
\begin{overpic}[scale=.5]{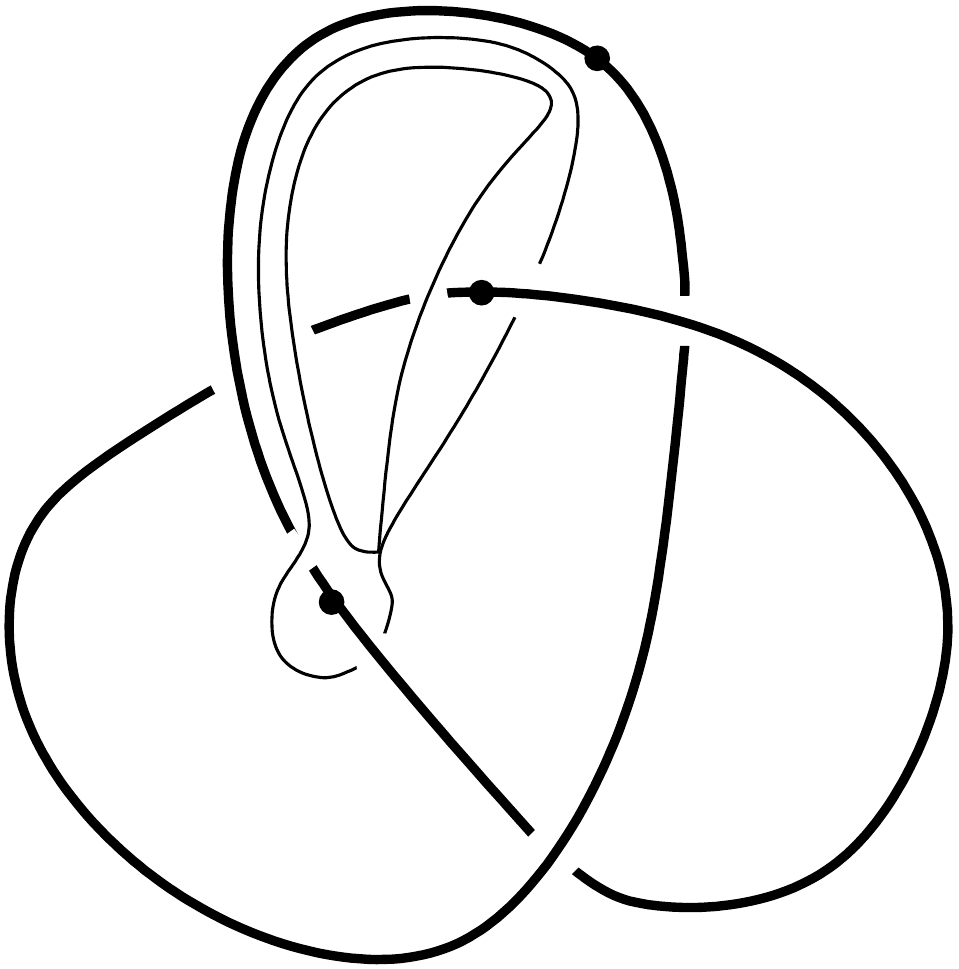}
\put(29.5,33){$a$}
\put(48.5,72.5){$b$}
\put(61,96.5){$c$}
\put(48,53){$\delta$}
\put(33,83){$\delta'$}
\end{overpic}\hfil
\begin{overpic}[scale=.5]{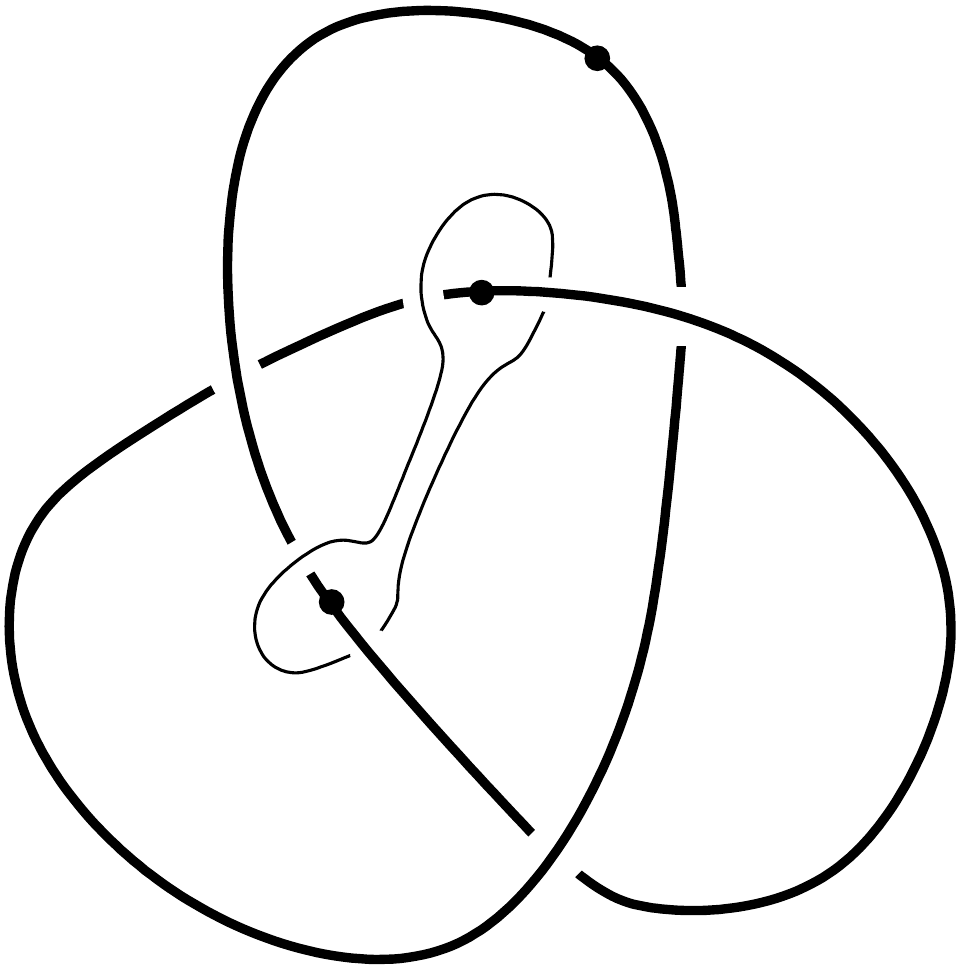}
\put(29.5,33){$a$}
\put(48.5,72.5){$b$}
\put(61,96.5){$c$}
\put(46,48){$l_{\segment{ab}}$}
\end{overpic}
}
\centerline{
\begin{overpic}[scale=.5]{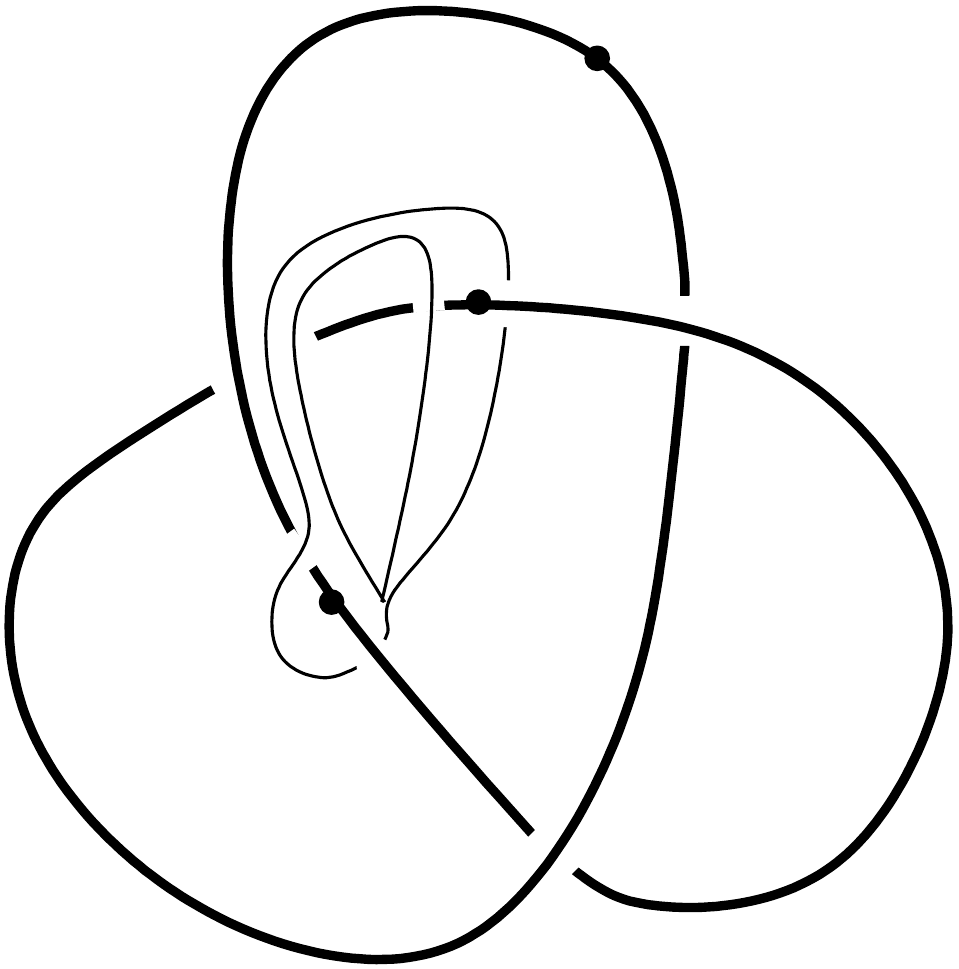}
\put(29.5,33){$a$}
\put(47.5,71.5){$b$}
\put(61,96.5){$c$}
\end{overpic}
}
\caption[The change from inessential to essential is a loop along a segment]
{At the top left, we see the loops $\delta$ and~$\delta'$, representing
the homotopy classes $h$ and~$h'$ arising from the essential
and inessential perturbations $S$ and~$S'$ 
of~$\segment{ac}$ as in \figr{secperturb}.
Because $\delta'$ is null-homotopic, applying the same motion to
parts of~$\delta$ shows it is homotopic 
to the loop~$l_{\segment{ab}}$ (top right) along~$\segment{ab}$.
At the bottom we see an intermediate stage of the homotopy.}
\label{fig:loop}
\end{figure}

In full generality, the secant $\segment{ac}$ may intersect
$K$ in many points (even infinitely many).  But still, the
two fundamental group elements $h$ and~$h'$ differ by some
finite word: the difference between~$S$ and~$S'$ is captured
by wrapping a different number of times around~$K$ at some finite number of 
intersection points $b_1, \ldots, b_k$, as in \figr{change}.
In particular, for some integers~$n_i$, we have that $S$ wraps
$n_i$ more times around~$b_i$ than~$S'$ does.  We can change from~$S$
to~$S'$ in $\sum n_i$ steps, at each step making just the simple
kind of change shown in \figr{secperturb}.
At (at least) one of these steps, we see a change from essential to inessential.
As in the simple case above, the homotopy class~$h$ just before such
a step can be represented by a loop~$l_{\segment{ab_i}}$ along
some segment~$\segment{ab_i}$.

\begin{figure}[ht!]\centering
\begin{overpic}[scale=.4]{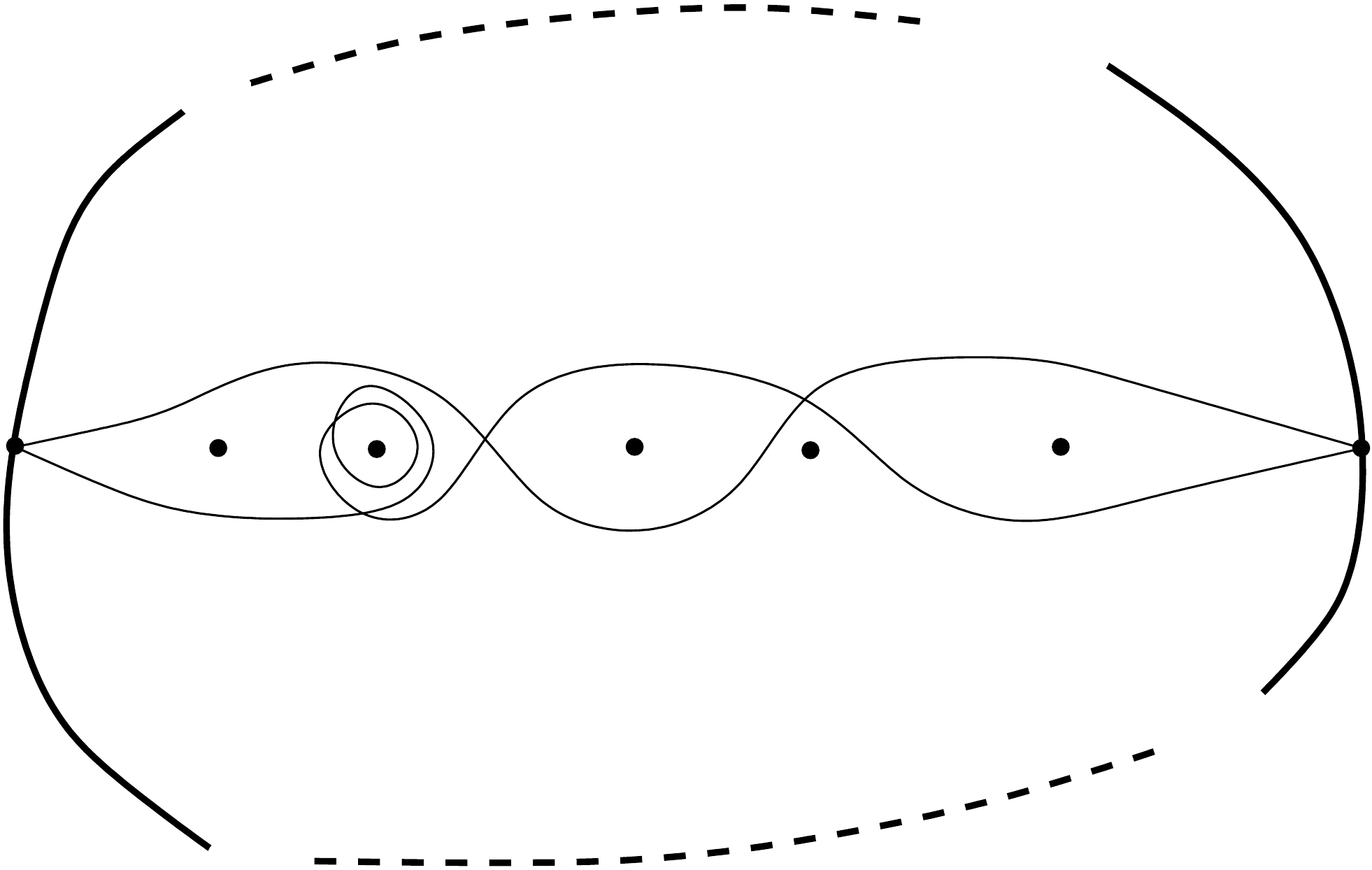}
{\small
\put(-3,30){$a$}
\put(15,38){$b_1$}
\put(27,37.5){$b_2$}
\put(45,38){$b_3$}
\put(58,38){$\cdots$}
\put(76,38){$b_k$}
\put(100.5,30){$c$}
\put(45,20){$S$}
\put(75,20){$S'$}
}
\end{overpic}
 \caption[The case of many intersection points]
{In case $\segment{ac}$ intersects~$K$ at many points, it still
has essential and inessential perturbations $S$ and $S'$,
and these differ by finitely many loops around intersection points $b_i$.}
\label{fig:change}
\end{figure}
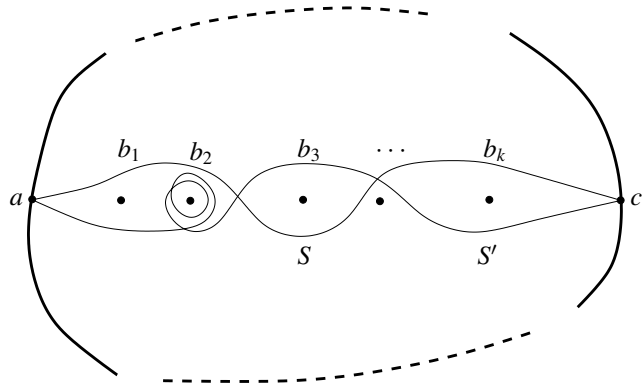

Note that because of the intersections (including $b_1,\ldots,b_{i-1}$)
of this segment with~$K$, the loop~$l_{\segment{ab_i}}$ is not \emph{a priori}
uniquely defined; it should be interpreted as wrapping around those previous
intersection points the same way the current $S$ does.  With this convention,
however, we see that $l_{\segment{ab_i}}$ represents the nontrivial homotopy
class~$h$.  The definition of strongly essential allows arbitrary small
perturbations, so again it follows immediately that $\segment{ab_i}$ is
strongly essential.  By symmetry, $\segment{b_ic}$ is also strongly essential.
Thus $b:=b_i$ is the desired intersection point for which $abc$ is essential.
\end{proof}

\section{Minimum arclength for essential subarcs of a knot}
We will improve our previous ropelength bounds by getting bounds
on the length of an essential arc.  A first bound is very easy:

\begin{lemma}\label{lem:a-b1}
If secant $ab$ is essential in a knot of unit thickness then $|a-b|\ge1$,
and if arc $\arc{ab}$ is essential then $\len{ab}\ge\pi$.
\end{lemma}

\begin{proof}
If $|a-b|<1$ then by \lem{thick} the ball $B$ of diameter $\segment{ab}$
contains a single unknotted arc (say $\arc{ab}$) of~$K$.
Now for any perturbation~$S$ of~$\segment{ab}$
which is disjoint from~$\arc{ab}$,
we can span $\arc{ab}\cup S$ by an embedded disk within $B$, whose
interior is then disjoint from~$K$.  This means that $\arc{ab}$ (and
thus $ab$) is inessential.

Knowing that sufficiently short arcs starting at any given point~$a$
are inessential, consider now the shortest arc $\arc{aq}$ which is
essential.  From \thm{change} there must be a trisecant $apq$
with both secants $ap$ and $pq$ essential, implying by the first
part that $a$ and $q$ are outside $B(p)$.  Since $ap$ is essential,
by the definition of~$q$ we have $p\notin\arc{aq}$,
meaning that $apq$ is reversed.
From \cor{diao} we get $\len{ab}\ge\len{aq}\ge\pi$.
\end{proof}

Intuitively, we expect an essential arc~$\arc{ab}$ of a knot to ``wrap
at least halfway around'' some point on the complementary arc~$\arc{ba}$.
Although when $|a-b|=2$ we can have~$\len{ab}=\pi$,
when $|a-b|<2$ we expect a better lower bound for $\len{ab}$.
Even though in fact an essential $\arc{ab}$ might instead ``wrap around''
some point on itself, we can still derive the desired bound.

\begin{lemma}\label{lem:lengthess}
If $\arc{ab}$ is an essential arc in a unit-thickness knot
and $|a-b|<2$, then~$\len{ab}\ge 2\pi-2\arcsin(|a-b|/2)$.
\end{lemma}

\begin{figure}[ht!]\centering
\begin{overpic}[scale=.5]{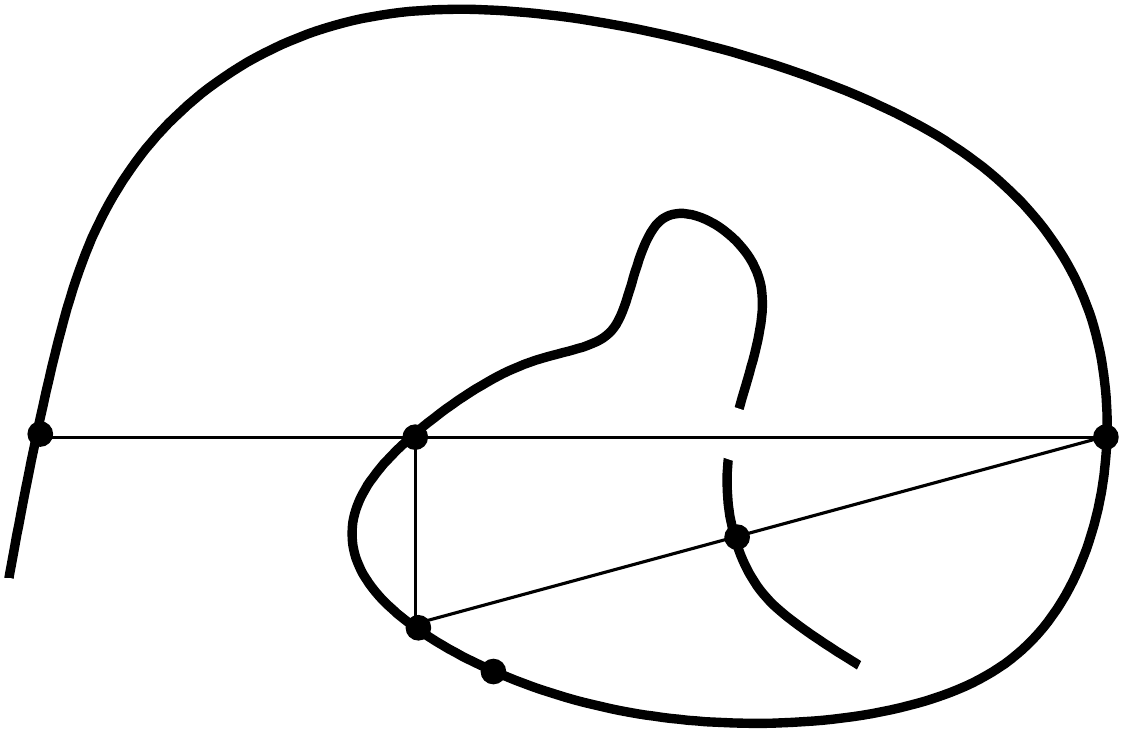}
{\small
\put(-2,25){$a$}
\put(35,30){$p$}
\put(101,25){$q$}
\put(38.5,16){$h$}
\put(60.5,12.5){$x$}
\put(32.3,4.8){$y$}
\put(42,-1){$b$}
}
\end{overpic}
\caption[Two essential trisecants]
{In the most intricate case in the proof of \lem{lengthess}, we let
$\arc{aq}$ be the first essential arc from~$a$, giving an essential
trisecant $apq$.  We then let $\arc{qy}$ be the first essential arc
from $q$, giving an essential trisecant $qxy$.  Since $|x-y|\ge 1$
and $|x-q|\ge 1$, setting $h=|p-y|$ we have $|p-q|\ge 2-h$.}
\label{fig:esslen}
\end{figure}
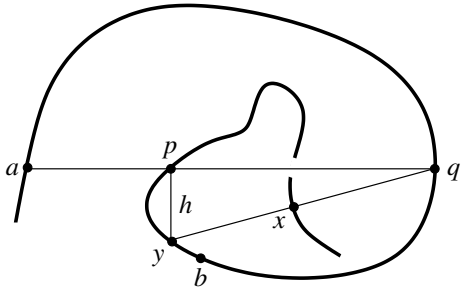

\begin{proof}
Note that $|a-b|\in[1,2]$, so $2\pi-2\arcsin(|a-b|/2)\le5\pi/3$.
As in the previous proof, let $\arc{aq}$ be the shortest
essential arc from~$a$, and find a reversed trisecant $apq$.
We have $b\notin\arc{aq}$ and $\len{aq}\ge\pi$,
so we may assume $\len{qb}<2\pi/3$ or the bound is trivially satisfied.

Since $\arc{qp}$ is essential, $\len{qp}\ge\pi>\len{qb}$, so
$b\in\arc{qp}$.  If $b\notin B(p)$
then the whole arc~$\arc{aqb}$ stays outside $B(p)$.
\begin{figure}[ht!]\centering
\begin{overpic}[scale=.4]{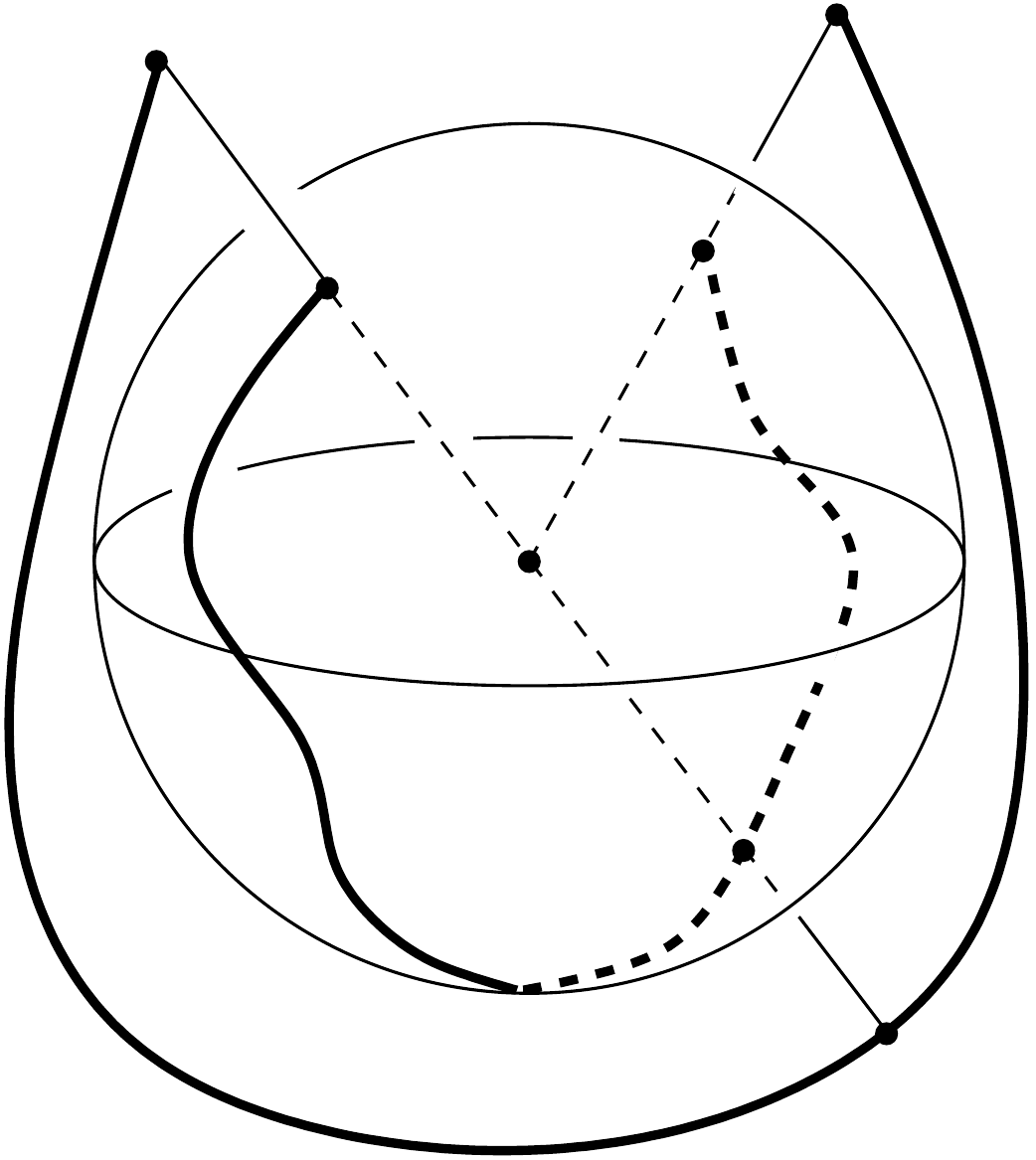}
{\small
\put(7.5,95){$a$}
\put(66.8,98){$b$}
\put(40.5,47.2){$p$}
\put(30,75){$\Pi a$}
\put(48.5,77.5){$\Pi b$}
\put(51.5,25.5){$\Pi q$}
\put(76,6){$q$}
}
\end{overpic}
\caption[Projection to the unit ball]
{In the proof of \lem{lengthess}, projecting $\arc{ab}$ to the unit ball
around~$p$ increases neither its length nor the distance between
its endpoints.  The projected curve includes antipodal points $\Pi a$
and $\Pi q$, which bounds its length from below.}
\label{fig:miness}
\end{figure}
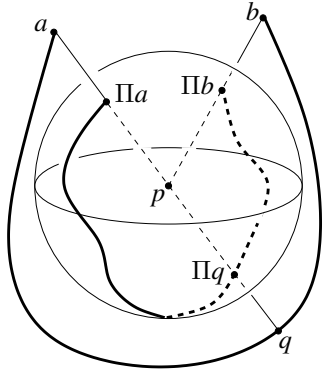
Let $\Pi$ denote the radial projection to
$\bdry B(p)$ as in \figr{miness}.
From \lem{balldist}, this projection does not increase length.
Because $|\Pi a-\Pi b|\le |a-b|$, we have
$2\pi-2\arcsin(|\Pi a-\Pi b|/2)\ge 2\pi-2\arcsin(|a-b|/2)$.
It therefore suffices to consider the case $\arc{ab}\subset \bdry B(p)$.
For any two points $x,y\in \bdry B(p)$, the spherical distance between them is
$2\arcsin(|x-y|/2)$. Thus
\begin{align*}
\len{ab} = \len{aq}+\len{qb} &\ge \pi+2\arcsin(|q-b|/2)\\
&= \pi+2\arccos(|a-b|/2) =2\pi-2\arcsin(|a-b|/2). 
\end{align*}

So we now assume that $|b-p|<1$.
Let~$\arc{qy}$ be the shortest essential arc starting from~$q$,
and note $|q-y|\ge 2$.
Since $\len{qy}\ge\pi>\len{qb}$ we have $b\in\arc{qy}$.
Let $h:=|p-y|\le\len{yp}$ and note that $h\in[0,1]$ since $b\in B(p)$.
(See \figr{esslen}.)
Since $|q-y|\ge2$, we have $|p-q|\ge2-h$, so $\len{aq}\ge \pi/2+f(2-h)$
by \cor{mindist}.
On the other hand, since $\len{bp}\le\pi/2$ (by \lem{thick})
and $\len{qy}\ge\pi$, we have
$\len{qb}\ge\pi/2+\len{yp}\ge\pi/2+h$.  Thus $\len{ab}\ge\pi+f(2-h)+h$.
An elementary calculation shows that the right-hand side is an increasing
function of~$h\in[0,1]$, minimized at $h=0$, where its value
is~$\pi+f(2)=7\pi/6+\sqrt3>5\pi/3$.
That is, we have as desired
$$\len{ab}\ge5\pi/3\ge 2\pi-2\arcsin(|a-b|/2).\proved$$
\end{proof}

If we define the continuous function
$$g(r) := \begin{cases}
              2\pi - 2\arcsin(r/2) &\text{if $0\le r\le 2$,} \\
              \pi            &\text{if $r\ge 2$.}
\end{cases}$$
then we can collect the results of the previous two lemmas as:

\begin{corollary}\label{cor:ess}
If $\arc{ab}$ is an essential arc in a knot~$K$ of unit thickness,
then
$$\len{ab} \ge g(|a-b|).$$
\end{corollary}

\section{Main results}
We now prove ropelength bounds for knots with different types
of quadrisecants. The following lemma will be used repeatedly.

\begin{lemma}\label{lem:mincalc}
Recall that
\begin{eqnarray*}
f(r)&:=&\sqrt{r^2-1}+\arcsin(1/r), \\
g(r)&:=& \begin{cases} 2\pi-2\arcsin(r/2) & \text{for } r\le2, \\
                        \pi & \text{for } r\ge2. \end{cases}
\end{eqnarray*}
Then, for $r\ge1$,
\begin{enumerate}
\item the minimum of $f(r)$ is $\pi/2$ and occurs at $r=1$,
\item the minimum of $f(r)+g(r)$ is $7\pi /6 + \sqrt{3}>5.397$ and occurs
  at $r=2$,
\item the minimum of $g(r)+r$ is $\pi+2>5.141$ and also occurs at $r=2$, and
\item the minimum of $2f(r)+g(r)+r$ is just over $9.3774$
  and occurs for $r\approx1.00305$.
\end{enumerate}
\end{lemma}

\begin{figure}[ht!]\centering
\centerline{
\begin{overpic}[scale=.7]{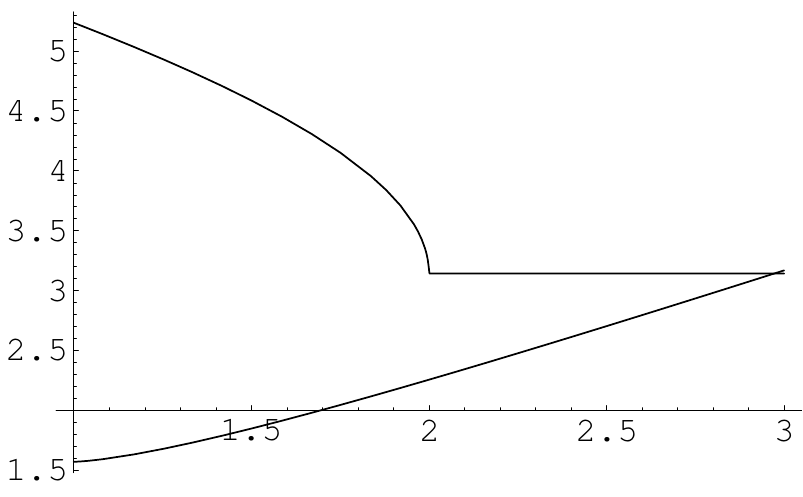}
\end{overpic}
\hfil
\begin{overpic}[scale=.7]{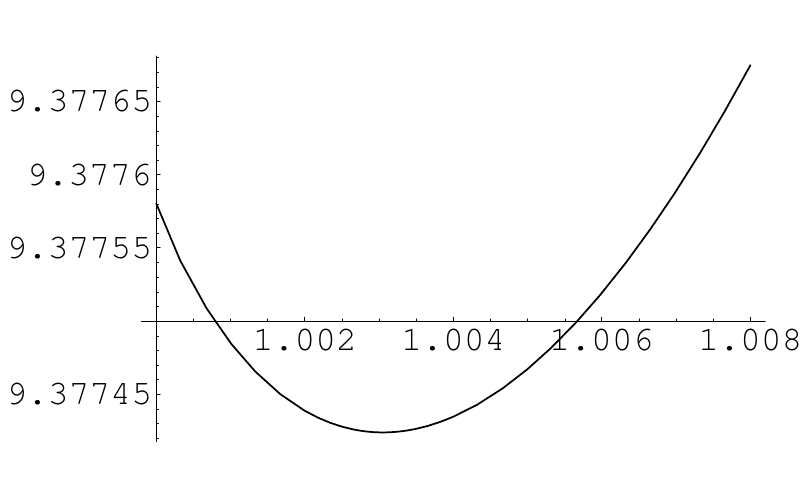}
\end{overpic}
}
\caption[Plots of the functions we minimize]
{Left, a plot of $f(r)$ and $g(r)$ for $r\in[1,3]$,
and right, a plot of $2f(r)+g(r)+r$ for $r\in[1,1.008]$.}
\label{fig:fgplot}
\end{figure}

\begin{proof}
  Note that $f$ is increasing, and $g$ is constant for $r\ge 2$.
  Thus the minima will occur in the range $r\in [1,2]$,
  where $f'=\tfrac1r\sqrt{r^2-1}$ and $g'=-2/\sqrt{4-r^2}$.
  Elementary calculations then give the results we want, where
  $r\approx1.00305$ is a polynomial root expressible in radicals.
  See also \figr{fgplot}.
\end{proof}

\begin{theorem}\label{thm:simple}
A knot with an essential simple quadrisecant has ropelength at least
$10\pi /3 +2\sqrt{3} +2> 15.936$.
\end{theorem}

\begin{proof}
Rescale the knot $K$ to have unit thickness, let $abcd$ be the
quadrisecant and orient $K$ in the usual way.
Then the length of~$K$ is
$\len{ab}+\len{bc}+\len{cd}+\len{da}$.
As before, let $r=|a-b|$, $s=|b-c|$ and $t=|c-d|$.

\fullref{cor:ess}
bounds $\arc{ab}$, $\arc{bc}$ and $\arc{cd}$.
The quadrisecant is essential, so from \lem{a-b1} we have~$r,s,t\ge 1$, and \lem{longarc} may be applied to bound $\len{da}$.
Thus the length of $K$ is at least
\begin{align*}
& g(r)+ g(s) + g(t)+ \big(f(r)+s + f(t)\big) \\
&\qquad=\quad \big(g(r)+f(r)\big) + \big(g(s)+s \big) + \big(g(t)+f(t)\big). 
\end{align*}

Since this is a sum of functions in the individual variables,
we can minimize each term separately.  These are the functions
considered in \lem{mincalc}, so the minima are achieved at $r=s=t=2$.
Adding the three values together, we find the ropelength of~$K$ is at least
$10\pi /3 +2\sqrt{3} +2 > 15.936$.
\end{proof}

\begin{theorem}\label{thm:flipped}
A knot with an essential flipped quadrisecant has ropelength at least
$10\pi/3 +2\sqrt{3}>\minval$.
\end{theorem}

\begin{proof}
Rescale $K$ to have unit thickness, let $abcd$ be the quadrisecant.
With the usual orientation, the length of~$K$ is
$\len{ab}+\len{bd}+\len{dc}+\len{ca}$.
Since the quadrisecant is essential,
from \lem{a-b1} and \lem{flipseg} we have $r,s,t\ge 1$.
We apply \cor{ess} to $\arc{ab}$ and
$\arc{dc}$ and \cor{mindist} to $\arc{bd}$ and $\arc{ca}$.

Thus the length of~$K$ is at least
\begin{align*}
& g(r) + \big(f(r)+f(s)\big) + g(t) + \big(f(s)+f(t)\big) \\
&\qquad=\quad \big(g(r) + f(r)\big) + 2f(s) + \big(g(t) + f(t)\big).
\end{align*}
Again we minimize the terms separately using \lem{mincalc}.
We find the ropelength of~$K$ is at least $10\pi /3 +2\sqrt{3}> \minval$.
\end{proof}

\begin{theorem}\label{thm:alternating}
A knot with an essential alternating quadrisecant has ropelength
at least~$\altval$.
\end{theorem}

\begin{proof}
Rescale $K$ to have unit thickness, let $abcd$ be the quadrisecant
and orient~$K$ in the usual way.
Then the ropelength of $K$ is $\len{ac}+\len{cb}+\len{bd}+\len{da}$.
Again, let $r=|a-b|$, $s=|b-c|$ and $t=|c-d|$.

The quadrisecant is essential, so from \lem{a-b1} and \lem{altarc}
we see $r,s,t\ge 1$.
Thus \lem{longarc} may be applied to $\arc{da}$.
We apply \cor{mindist} to $\arc{ac}$ and $\arc{bd}$,
and \cor{ess} to $\arc{cb}$.

We find that the length of~$K$ is at least
\begin{align*}
&\big(f(r)+f(s)\big) + \big(f(s)+f(t)\big)
    + g(s) + \big(f(r)+s+f(t)\big) \\
&\qquad =\quad 2f(r) + \big(2f(s) + g(s) + s\big) + 2f(t).
\end{align*}

Again, we can minimize in each variable separately, using \lem{mincalc}.
Hence the ropelength of~$K$ is at least $2\pi + 9.377 > \altval$.
\end{proof}

\begin{theorem}
Any nontrivial knot has ropelength at least $\minval$.
\end{theorem}

\begin{proof}
Any knot of finite ropelength is $C^{1,1}$, so by \cor{alless}
it has an essential quadrisecant.
This must be either simple, alternating, or flipped, so
one of the theorems above applies; we inherit the
worst of the three bounds.
\end{proof}

In her doctoral dissertation~\cite{denne}, Denne shows:
\begin{theorem}\label{thm:essaltquad}
Any nontrivial $C^{1,1}$ knot has an essential alternating quadrisecant.
\end{theorem}

Combining this with \thm{alternating} gives:
\begin{corollary}
Any nontrivial knot has ropelength at least $\altval$.
\end{corollary}

We note that this bound is better even than the conjectured
bound of 15.25 from
\hbox{\cite[Conjecture~26]{cks}}.
We also note that our bound cannot be sharp, for a curve
which is~$C^1$ at~$b$ cannot simultaneously achieve the
bounds for~$\len{cb}$ and~$\len{bd}$ when $s\approx1.003$.
Probably a careful analysis based on the tangent directions
at~$b$ and~$c$ could yield a slightly better bound.
However, we note again that numerical simulations have found trefoil
knots with ropelength no more than 5\% greater than our bound,
so there is not much further room for improvement.

\section{Quadrisecants and links}

Quadrisecants may be used in a similar fashion to give
lower bounds for the ropelength of a nontrivial link.
For a link, $n$--secant lines and $n$--secants are defined as before.
They can be classified in terms of the order in which they intersect
the different components of the link.
We begin our considerations with a simple lemma bounding
the length of any nonsplit component of a link.

\begin{lemma}
Suppose $L$ is a link, and $A$ is a component of~$L$ not split
from the rest of the link.  Then for any point $a\in A$,
there is a trisecant $aba'$
where $a'\in A$ but $b$ lies on some other component.
If~$L$ has unit thickness, then $A$ has length at least~$2\pi$.
\end{lemma}
\begin{proof}
For the given~$a$, the union of all secants~$aa'$ is a disk
spanning~$A$.  Because~$A$ is not split from~$L\setm A$,
this disk must be cut by $L\setm A$, at some point~$b$.
This gives the desired trisecant.  If $L$ has unit thickness,
then $A$ stays outside $B(b)$, so as in \cor{diao}, we have
$\len{aa'}\ge\pi$ and $\len{a'a}\ge\pi$.
\end{proof}

This construction of a trisecant is adapted from
Ortel's original solution of the Gehring link problem.  (See~\cite{cfksw}.)
The length bound can also be viewed as a special case of
\hbox{\cite[Theorem~10]{cks}},
and immediately implies the following corollary.

\begin{corollary}\label{cor:nonsplit}
A nonsplit link with $k$ components has ropelength at least~$2\pi k$.
\end{corollary}

This bound is sharp in the case of the tight Hopf link, where each
component has length exactly~$2\pi$.

We know that nontrivial links have many trisecants, and want to consider when
they have quadrisecants.
Pannwitz~\cite{pann} was the first to show the existence of
quadrisecants for certain links:

\begin{proposition} 
If $A$ and $B$ are disjoint generic polygonal knots, linked in the sense that
neither one is homotopically trivial in the complement of the other,
then $A\cup B$ has a quadrisecant line intersecting them in the order $ABAB$.
\end{proposition}

Note that, when $A$ or $B$ is unknotted, Pannwitz's hypothesis
is equivalent to having nonzero linking number.  The theorem says
nothing, for instance, about the Whitehead link.

Kuperberg~\cite{kup} extended this result to apply
to all nontrivial link types (although with no information
about the order type), and for generic links
he again guaranteed an essential quadrisecant:

\begin{proposition}
A generic nontrivial link has an essential quadrisecant.
\end{proposition}

Here, a secant of a link~$L$ is automatically essential if
its endpoints lie on different components of~$L$.  If its
endpoints are on the same component~$K$, we apply our
previous definition of inessential secants, but require the disk
to avoid all of~$L$, not just the component~$K$.

Combining Kuperberg's result with \prop{limess}
(which extends easily to links) immediately gives us:

\begin{theorem}\label{thm:linkquad}
Every nontrivial $C^{1,1}$ link has an essential quadrisecant.
\end{theorem}

Depending on the order in which the quadrisecant visits the
different components of the link, we can hope to get bounds
on the ropelength.  Sometimes, as for $ABAB$ quadrisecants,
we cannot improve the bound of \cor{nonsplit}.  Here of course,
the example of the tight Hopf link (which does have an $ABAB$
quadrisecant) shows that this bound has no room for improvement.

Note that \lem{a-b1} extends immediately to links, since if
$a$ and~$b$ are on different components of a unit-thickness link
we automatically have $|a-b|\ge1$.

\begin{theorem}\label{thm:linklen}
Let $L$ be a link of unit thickness with an essential quadrisecant~$Q$.
If~$Q$ has type $AAAB$, $AABA$, $ABBA$ or $ABCA$,
then the length of component~$A$ is at least
$7\pi/3+2\sqrt{3}$, $8\pi/3+1+\sqrt{3}$, $2\pi+2$ or $2\pi+2$, respectively.
\end{theorem}

\begin{proof}
Suppose $a_1a_2a_3b$ is an essential quadrisecant of type $AAAB$,
and set
$$r:=|a_2-a_1|, \qquad s:=|a_3-a_2|.$$
The quadrisecant is
essential, so from \lem{a-b1} we have $r,s\ge 1$ and we may
apply \cor{mindist} to $\arc{a_1a_3}$ and \cor{ess} to $\arc{a_1a_2}$
and $\arc{a_2a_3}$. Thus the length of~$A$ is at least $f(r)+g(r) + f(s)+g(s)$.
As before, we can minimize in each variable separately, using \lem{mincalc}.
Hence the length of~$A$ is at least $7\pi/3 +2\sqrt{3}$.
 
Now suppose that $a_1a_2ba_3$ is an essential quadrisecant of type $AABA$,
and set
$$r:=|a_1-a_2|, \qquad s:=|a_2-b|, \qquad t:=|b-a_3|.$$
Since the
quadrisecant is essential, we have $r,s,t\ge 1$ and we may
apply \cor{ess} to $\arc{a_1a_2}$, \cor{mindist} to $\arc{a_2a_3}$ and
\lem{longarc} to $\arc{a_1a_3}$. 
We find that the length of $A$ is at least
\begin{align*}
&g(r) +\big(f(s)+f(t)\big)+ \big(f(r)+s + f(t)\big) \\
&\qquad =\quad \big(f(r)+g(r)\big) + \big(f(s) + s\big) + 2f(t).
\end{align*}
Again, minimizing in each variable separately using \lem{mincalc},
we find the length of~$A$ is at least $8\pi/3+1+\sqrt{3}$.

Finally suppose $a_1bca_2$ is an essential quadrisecant of type $ABBA$
or type $ABCA$.  Because the quadrisecant is essential
we must have $|a_1-b|\ge 1$, $|b-c|\ge 1$
and $|c-a_2|\ge 1$. Just as in the proof of \lem{longarc},
we find that $\len{a_1a_2}\ge \pi+1$ and that~$\gamma_{a_2a_1}\ge \pi+1$,
showing that the length of~$A$ is at least $2\pi+2$ as desired.
\end{proof}

Unfortunately, we do not know any link classes which would be
guaranteed to have one of these types of quadrisecants, so we
know no way to apply this theorem.

\subsection*{Acknowledgements}
We extend our thanks to Stephanie Alexander, Dick Bishop,
Jason Cantarella, Rob Kusner and Nancy Wrinkle
for helpful conversations and suggestions.
Our work was partially supported by the National Science Foundation
through grants DMS--03--10562 (Diao) and DMS--00--71520 (Sullivan and Denne),
and by the UIUC Campus Research Board.


\bibliographystyle{gtart}
\bibliography{link}

\end{document}